\documentclass{amsproc}

\usepackage{graphics}

\newtheorem{theorem}{Theorem}[section]
\newtheorem{lemma}[theorem]{Lemma}
\newtheorem{corollary}[theorem]{Corollary}
\newtheorem{proposition}[theorem]{Proposition}
\newtheorem{definition}[theorem]{Definition}
\newtheorem{example}[theorem]{Example}
\newtheorem{remark}[theorem]{Remark}

\numberwithin{equation}{section}

\newcommand{\bz}{{\mathbb B}}
\newcommand{\cz}{{\mathbb C}}
\newcommand{\gz}{{\mathbb Z}}

\newcommand{\nz}{{\mathbb N}}
\newcommand{\rz}{{\mathbb R}}

\newcommand{\calA}{\mathcal{A}}
\newcommand{\calC}{\mathcal{C}}
\newcommand{\calD}{\mathcal{D}}
\newcommand{\calE}{\mathcal{E}}
\newcommand{\calF}{\mathcal{F}}
\newcommand{\calH}{\mathcal{H}}
\newcommand{\calK}{\mathcal{K}}
\newcommand{\calL}{\mathcal{L}}
\newcommand{\calM}{\mathcal{M}}

\newcommand{\calS}{\mathcal{S}}

\newcommand{\amin}{A_{\text{\rm min}}}
\newcommand{\bzbz}{{\mathbb B}\times{\mathbb B}}
\newcommand{\ci}{\mathcal{C}^\infty}
\newcommand{\cii}{\mathcal{C}^{\infty,\infty}}
\newcommand{\cicomp}{\mathcal{C}^\infty_{\text{\rm comp}}}

\newcommand{\dbar}{d\hspace*{-0.08em}\bar{}\hspace*{0.1em}}

\newcommand{\intb}{\text{\rm int}\,{\mathbb B}}
\newcommand{\op}{\text{\rm op}}

\newcommand{\pit}{\,{\widehat{\otimes}}_\pi\,}
\newcommand{\re}{\text{\rm Re}\,}
\newcommand{\rpbar}{\overline{{\mathbb R}}_+}
\newcommand{\skp}[2]{\langle#1,#2\rangle}
\newcommand{\skpa}[2]{[#1,#2]}
\newcommand{\smsum}{\mathop{\mbox{$\sum$}}}
\newcommand{\spk}[1]{\left<#1\right>}
\newcommand{\st}{\mbox{\boldmath$\;|\;$\unboldmath}}
\newcommand{\trinorm}[1]%
    {|\hspace*{-1pt}|\hspace*{-1pt}|#1|\hspace*{-1pt}|\hspace*{-1pt}|}

\renewcommand{\Re}{{\rm Re}\,}

\begin{document}

\title[The Resolvent of Closed Extensions of Cone Differential Operators]
      {The Resolvent of Closed Extensions of Cone Differential Operators}
%%%%%%%%%%%%%%%%%%%%%%%%%%%%%%%%%%%%%%%%%%%%%%%%%%%%%%%%%%%%%%%%%%%%%
\author{E.\ Schrohe}
\address{Universit\"at Potsdam, Institut f\"ur Mathematik,
         Postfach 60 15 53, 14415 Potsdam, Germany}
\email{schrohe@math.uni-potsdam.de}
\author{J.\ Seiler}
\address{Universit\"at Potsdam, Institut f\"ur Mathematik,
         Postfach 60 15 53, 14415 Potsdam, Germany}
\email{seiler@math.uni-potsdam.de}
\date{\today}
%%%%%%%%%%%%%%%%%%%%%%%%%%%%%%%%%%%%%%%%%%%%%%%%%%%%%%%%%%%%%%%%%%%%%
%%%%%%%%%%%%%%%%%%%%%%%%%%%%%%%%%%%%%%%%%%%%%%%%%%%%%%%%%%%%%%%%%%%%%
\begin{abstract}
 We study closed extensions $\underline A$ of
 an elliptic differential operator $A$ on a manifold with conical
 singularities, acting as an unbounded operator on a weighted $L_p$-space.
 Under suitable conditions we show that the resolvent 
 $(\lambda-\underline A)^{-1}$ exists
 in a sector of the complex plane and decays like $1/|\lambda|$ as
 $|\lambda|\to\infty$. Moreover, we determine the structure of the resolvent
 with enough precision to guarantee existence and boundedness of imaginary
 powers of $\underline A$. 
 
 As an application we treat the Laplace-Beltrami operator for a metric with
 straight conical degeneracy and describe domains yielding
 maximal regularity for the Cauchy problem $\dot{u}-\Delta u=f$, $u(0)=0$. 
\end{abstract}

\maketitle

%\begin{center}
% Preliminary version: \quad \today
%\end{center}

%%%%%%%%%%%%%%%%%%%%%%%%%%%%%%%%%%%%%%%%%%%%%%%%%%%%%%%%%%%%%%%%%%%%%
%%
%% To insert footnote on first page but not on all subsequent ones
%%
%%\markboth{\uppercase{
%%          The resolvent of Extensions of Cone Differential Operators}}
%%         {\uppercase{E.\ Schrohe and J.\ Seiler}}
%%
%%%%%%%%%%%%%%%%%%%%%%%%%%%%%%%%%%%%%%%%%%%%%%%%%%%%%%%%%%%%%%%%%%%%%
\tableofcontents
%%%%%%%%%%%%%%%%%%%%%%%%%%%%%%%%%%%%%%%%%%%%%%%%%%%%%%%%%%%%%%%%%%%%%

\section{Introduction}\label{intro}

Understanding the resolvent of elliptic differential operators is of central
interest for many questions in partial differential equations. Following the
approach suggested by Seeley, it is crucial for the
analysis of the heat operator or of complex powers.
In his classical paper \cite{Seel0},
he showed how the parametrix to an elliptic operator on a closed manifold
can be constructed as a parameter-dependent pseudodifferential operator and how
the structure of the parametrix determines the
essential properties of the complex powers.
He subsequently extended his methods to cover
also boundary value problems \cite{Seel1}
and proved the boundedness of the purely
imaginary powers \cite{Seel2}. 
His results have attracted new interest
in connection with modern methods in nonlinear
evolution equations, where one requires maximal regularity for the
generator of the associated semigroup, 
which in turn is implied by the boundedness of its purely imaginary powers. 

In the present paper we study an elliptic differential operator $A$ on a
manifold $B$ with conical singularities (a `cone differential operator').
The investigation of these operators started with the work of Cheeger
\cite{Che}. Important contributions to the index theory were made 
in particular by
Br\"uning \& Seeley \cite{BrSe1} and Lesch \cite{Lesc}; 
associated pseudodifferential calculi were devised by
Melrose \cite{RBM}, Plamenevskij \cite{Plam86}, and Schulze \cite{Schu2}.

While the picture of the conical singularity helps the intuition, one 
prefers to perform the actual analysis on a manifold $\bz$ with boundary,
thought of as the blow-up of $B$. 
A cone differential operator of order $\mu$ by definition is an
operator that can be written  in the form
$$A = t^{-\mu}\smsum_{j=0}^\mu a_j(t)(-t\partial_t)^j$$
in a neighborhood of the boundary. 
Here $t$ is a boundary defining function 
and $a_j$ a smooth family of differential operators of order $\mu-j$ 
on $\partial\bz$.

We consider $A$ as an unbounded operator acting in a (weighted) $L_p$-space. 
Our goal is to find conditions which ensure the existence of the resolvent $(A-
\lambda)^{-1}$ in a sector of the complex plane with decay like $1/|\lambda|$
as $|\lambda|\to\infty$ and to determine its structure with enough precision to
construct complex powers and to show their boundedness for purely imginary
exponents. We work with a variant of Schulze's cone calculus
because the concept of meromorphic Mellin symbols makes it easy to
describe the connection between operators and function spaces with asymptotics.

A cone differential operator in general has
many closed extensions, see e.g.\ Lesch \cite[Section 1.3]{Lesc}.
While, {\em a priori}, there is no preference for any of these from
the analytical point of view, it is obvious that only for few of them the
resolvent will have good properties. One basic problem therefore is to 
determine all possible choices.
Our Theorem \ref{maxdomain} completes Lesch's results in that we obtain an
explicit formula for the domain of the maximal extension in the general
situation. 

Extending Theorem 3.14 from \cite{ScSe}, we next clarify the structure of the
inverse of a bijective closed extension of $A$. In  Theorem \ref{propb} we show
how $A^{-1}$ can be decomposed as the sum of two operators in usual cone
calculi for different weight data.

We then turn to the analysis of the resolvent. In order to keep the exposition 
short and the proofs transparent we restrict ourselves to the case where 
the coefficients $a_j$ of the operator $A$ are constant for small $t$.
The general case will be treated in a subsequent publication.

Following a standard technique, we replace the spectral parameter 
$\lambda$ by $\eta^\mu$, where $\mu$ is the order of $A$, and $\eta$
varies in a corresponding sector of $\cz$.    
In close analogy to Theorem \ref{propb} we prove in Theorem \ref{maintheorem}
that $(A-\eta^\mu)^{-1}$ is 
the sum of two parameter-dependent cone operators; the parameter space is
the new $\eta$-sector.
In order to establish this fact we have to make assumptions
which are restrictive but nevertheless seem natural in this context:
Clearly, we have to ask for the invertibility of the principal 
pseudodifferential symbol of $A-\eta^\mu$ in the sector, 
with a certain uniformity as one approaches the singularity. 
Moreover, we require the invertibility of  $\widehat A- \eta^\mu$, where 
$\widehat A$ is the `model cone operator' associated to $A$. 
It is given by $\widehat{A}=t^{-\mu}\sum_{j=0}^\mu
a_j(0)(-t\partial_t)^j$ on $\rz_+\times\partial\bz$ and reflects 
the behavior of $A$ near the singular point; $\widehat{A}$ acts on a domain 
linked to that of $A$. As $\widehat A-\eta^\mu$ can be considered the analog of
an edge principal symbol for $A-\eta^\mu$, its invertibility appears to be
necessary for the above result. Finally, we assume for technical reasons that
the domain of $A$ (or more precisely the associated domain of  $\widehat A$) 
is invariant under dilations (`saturated' in the language of Gil and Mendoza
\cite{GiMe}).  

It follows from Theorem 5.1 and Remark 5.5 in \cite{CSS1}
that the structure of 
the resolvent we obtain from Theorem \ref{maintheorem} is
precisely that required for the construction of complex powers and implies
the boundedness of the purely imaginary powers; we can hence extend the 
results of that paper as well as those in \cite{CSS2} to this new class of 
operators.

The idea of analyzing the resolvent of a cone differential operator
in terms of a suitable pseudodifferential calculus is not new.
In fact, writing the resolvent as a parameter-dependent 
cone operator can be seen as a special case of the edge parametrix 
construction, see Schulze \cite[Section 9.3.3, Theorem 6]{EgSc}.
Moreover, Gil \cite{Gil}, \cite{GilMN},
and Loya \cite{Loya1}, \cite{Loya2}, also 
in joint work \cite{GilLoya}, used this technique to derive results on
heat invariants, complex powers, and noncommutative residues.
While these are important theorems, there is one draw-back:
In all articles, the authors rely on a special form
of the above ellipticity condition, namely the
invertibility of $\widehat A-\eta^\mu$, acting between weighted Mellin Sobolev 
spaces. 
One can show, however, that this assumption fails in many cases, e.g.\ 
for the Laplace-Beltrami operator in dimensions $\le4$, acting in $L^2$
with respect to any metric that has a straight conical singularity.
Roughly speaking, this approach works only for the minimal (and hence by 
duality for the maximal) extension.  
The new point here is that we can now treat all closed extensions with dilation 
invariant (saturated) domains, opening the way for the analysis of larger
classes of operators.  
 
%Mooers in \cite{Moo} investigated the structure of the 
%heat kernel of self-adjoint extensions of 
%the Laplacian on forms.
%In the context of heat kernel expansions Mooers in \cite{Moo} 

As an application we study the Laplacian in weighted $L_p$-spaces,
$1<p<\infty$. Combining our analysis with techniques of Gil and Mendoza
\cite{GiMe}, we show in Theorems \ref{elldomain1} and \ref{elldomain2} how one
can always choose the domain in such a way that the above ellipticity
conditions are fulfilled. This yields maximal regularity for the Cauchy 
problem $\dot{u}-\Delta u=f$ {on} $] 0,T[$, $u(0)=0,$ which is the starting
point for many results in nonlinear evolution equations.

%%%%%%%%%%%%%%%%%%%%%%%%%%%%%%%%%%%%%%%%%%%%%%%%%%%%%%%%%%%%%%%%%%%%%%
%%%%%%%%%%%%%%%%%%%%%%%%%%%%%%%%%%%%%%%%%%%%%%%%%%%%%%%%%%%%%%%%%%%%%%

\section{Cone differential operators and their closed extensions}

%%%%%%%%%%%%%%%%%%%%%%%%%%%%%%%%%%%%%%%%%%%%%%%%%%%%%%%%%%%%%%%%%%%%%%

\subsection{Operators on $\bz$}\label{section1.1}
 
 Let $\bz$ be a smooth, compact manifold with boundary. A $\mu$-th order
 differential operator $A$ with smooth coefficients acting on sections of a
 vector bundle $E$ over the interior of $\bz$ is called a {\em cone
 differential operator} if, near the boundary, it has the form
  \begin{equation}\label{coneoperator}  
   A=t^{-\mu}\smsum_{j=0}^\mu a_j(t)(-t\partial_t)^j,\qquad 
   a_j\in\ci([0,1[,\mbox{\rm Diff}^{\mu-j}(\partial\bz)).
  \end{equation}
 In more detail: We assume (as we may) that $E$ respects the product structure
 near the boundary, i.e.\ $E$ is the pull-back of a vector bundle $E_\partial$
 over $\partial\bz$ under the canonical projection
 ${[0,1[}\times\partial\bz\to\partial\bz$. The coefficients $a_j(t)$ then are
 differential operators acting on sections of $E_\partial$. In order to keep
 the exposition simple, however, we shall not indicate the bundles in the
 notation. 
 
 Besides the standard pseudodifferential principal symbol
 $\sigma^\mu_\psi(A)\in\ci(T^*\intb\setminus0)$, we associate with $A$ 
 two other symbols: First, there is the {\em rescaled symbol}
 $\widetilde{\sigma}^\mu_\psi(A)\in
 \ci((T^*\partial\bz\times\rz)\setminus0)$ which, in local coordinates, is
 given by 
  $$\widetilde{\sigma}^\mu_\psi(A)(x,\xi,\tau)=\smsum_{j=0}^\mu 
    \sigma^{\mu-j}_\psi(a_j)(0,x,\xi)(-i\tau)^j.$$
 Secondly, we have the 
 {\em conormal symbol} $\sigma^\mu_M(A)$ defined by
  $$\sigma^\mu_M(A)(z)=\smsum_{j=0}^\mu a_j(0)z^j,\qquad z\in\cz.$$
 It is a polynomial in $z$ of degree at most $\mu$ with values in differential
 operators on $\partial\bz$ of order at most $\mu$. In particular, 
 $\sigma^\mu_M(A)\in
 \calA(\cz,\calL(H^s_p(\partial\bz),\calL(H^{s-\mu}_p(\partial\bz)))$ for all
 $s\in\rz$, $1<p<\infty$, where $\calA(\cz,X)$ denotes the holomorphic,
 $X$-valued functions on $\cz$. 
 
 Let us introduce some notions we shall frequently use throughout this paper. 
 
 A {\em cut-off function}, generally denoted by $\omega$, $\omega_0$,
 $\omega_1$, or $\sigma$, $\sigma_0$, $\sigma_1$, is a non-negative decreasing
 function in $\cicomp([0,1[)$, which is identically 1 near zero. 
 
 \begin{definition}\label{belliptic}
  \begin{itemize}
   \item[a)] $A$ is called $\bz$-elliptic if both   
    $\sigma^\mu_\psi(A)$ and $\widetilde{\sigma}^\mu_\psi(A)$ are pointwise
    invertible. 
   \item[b)] $A$ is said to have $t$-independent coefficients near
    the boundary if the functions $a_j$ in \eqref{coneoperator} are constant
    in $t$. 
  \end{itemize}
 \end{definition}
 
 The operator $A$ induces continuous actions 
  \begin{equation}\label{bounded}
   A:\calH^{s,\gamma}_p(\bz)\longrightarrow\calH^{s-\mu,\gamma-\mu}_p(\bz),
   \qquad
   s,\gamma\in\rz,\;1<p<\infty,
  \end{equation}
 in a scale of Sobolev spaces which is defined as follows: 
 
 \begin{definition}\label{sobolev}
  Let $s\in\nz_0$. The space of all distributions $u\in H^s_{p,loc}(\intb)$
  with
   $$t^{\frac{n+1}{2}-\gamma}(t\partial_t)^k\partial^\alpha_x(\omega u)(t,x)\in
     L_p([0,1[\times\partial\bz,\mbox{$\frac{dt}{t}dx$})\qquad
     \forall\;k+|\alpha|\le s$$ 
  is denoted by $\calH^{s,\gamma}_p(\bz)$. Here, $\omega$
  denotes an arbitrary cut-off function. 
 \end{definition}
 
 This definition extends to real $s$, yielding a scale of Banach
 spaces (Hilbert spaces in case $p=2$) with two properties we want to
 mention explicitly: The embedding $\calH^{s^\prime,\gamma^\prime}_p(\bz)
 \hookrightarrow\calH^{s,\gamma}_p(\bz)$ is continuous for 
 $s^\prime\ge s$, $\gamma^\prime\ge\gamma$, and compact if 
 $s^\prime>s$, $\gamma^\prime>\gamma$; the scalar-product of
 $\calH^{0,0}_2(\bz)$ induces an identification of the dual space
 $(\calH^{s,\gamma}_p(\bz))^\prime$ with $\calH^{-s,-\gamma}_{p^\prime}(\bz)$,
 where $p^\prime$ is the dual number to $p$, i.e.\
 $\frac{1}{p}+\frac{1}{p^\prime}=1$. 

 Instead of considering $A$ as a continuous operator in the Sobolev spaces,
 %cf.\ \eqref{bounded}, 
 we shall now study the closed extensions of the 
 unbounded operator 
  \begin{equation}\label{unbounded}
   A:\cicomp(\intb)\subset\calH^{0,\gamma}_p(\bz)\longrightarrow 
   \calH^{0,\gamma}_p(\bz).
  \end{equation}
 %(without any difference we could also consider $A$ on $\cii(\bz)$, cf.\
 %Section \ref{appendix1}). A natural assumption we will impose on $A$ in the
 In the sequel  $A$ will be assumed to be $\bz$-elliptic and of positive order $\mu>0$.
 In the upcoming Sections \ref{section1.2},
 \ref{section1.4} we shall give an explicit description of all possible closed
 extensions of $A$. We shall need a few basic facts about the cone calculus
 which may be found in the short introduction \cite{Seil2}. We refer in particular
 to \cite{Seil2}, Section 2.4,  for  the notion of 
 Mellin pseudodifferential operators $\op_M^{\delta}$ and their mapping properties on the 
 spaces $\calH^{s,\gamma}_p(\bz)$. Some material can also be found in the appendix; 
these two points, for example are covered by \eqref{wxyz} and Remark
\ref{mapprop}.
 We let $\omega$ denote an arbitratry cut-off  function. 

%%%%%%%%%%%%%%%%%%%%%%%%%%%%%%%%%%%%%%%%%%%%%%%%%%%%%%%%%%%%%%%%%%%%%

\subsection{The minimal extension}\label{section1.2}

 The following result was shown in \cite{GiMe}, Proposition 3.6. We give here a
 short proof, using some results of \cite{Lesc}. 

 \begin{proposition}\label{mindomain}
  The domain of the closure $A_{\min}=A^{\gamma,p}_{\min}$ of $A$ from
  \eqref{unbounded} is 
  \begin{align}
   \calD(A_{\min})&=\calD(A_{\max})\cap
   \mathop{\mbox{\Large$\cap$}}_{\varepsilon>0}\,
   \calH^{\mu,\gamma+\mu-\varepsilon}_p(\bz)\nonumber\\
   &=\Big\{u\in\mathop{\mbox{\Large$\cap$}}_{\varepsilon>0}\,
   \calH^{\mu,\gamma+\mu-\varepsilon}_p(\bz)\st 
   t^{-\mu}\smsum_{j=0}^\mu a_j(0)(-t\partial_t)^j(\omega u)\in
   \calH^{0,\gamma}_p(\bz)\Big\}.\label{md}
  \end{align}
  In particular, 
   $$\calH^{\mu,\gamma+\mu}_p(\bz)\hookrightarrow\calD(A_{\min})
     \hookrightarrow\calH^{\mu,\gamma+\mu-\varepsilon}_p(\bz)\qquad
     \forall\;\varepsilon>0.$$
  We have $\calD(A_{\min})=\calH^{\mu,\gamma+\mu}_p(\bz)$ if and only if the
  conormal symbol $\sigma^\mu_M(A)(z)$ is invertible for all $z$ with 
  $\re z=\frac{n+1}{2}-\gamma-\mu$. 
 \end{proposition}
 
 \begin{proof}
  According to \cite{Lesc}, Proposition 1.3.12, we may assume that $A$ has
  $t$-independent coefficients near the boundary. Now let
  $u\in\calD(A_{\min})$, i.e.\ there exists a sequence of functions
  $u_n\in\cicomp(\intb)$ such that $u_n\to u$ and $Au_n\to Au$ with convergence
  in $\calH^{0,\gamma}_p(\bz)$. Choose a cut-off function 
  $\widetilde{\omega}$ with
  $\omega\widetilde{\omega}=\omega$, and let 
  $B=\widetilde{\omega}\,\op_M^{\gamma+\mu-\varepsilon-\frac{n}{2}}
  (\sigma^\mu_M(A)^{-1})\,t^\mu$ with arbitrarily small $\varepsilon>0$ (and
  $\varepsilon=0$ in case of the invertibility of the conormal symbol). 
  It follows from elliptic regularity that $u_n\to u$ in 
  $H^\mu_{p,{\rm loc}}(\intb)$. Thus $(1-\omega)u_n\to(1-\omega)u$ in
  $\calH^{\mu,\gamma+\mu}_p(\bz)$ and $A(\omega u_n)\to A(\omega u)$ in
  $\calH^{0,\gamma}_p(\bz)$. Therefore $\omega u$ belongs to $\calD(A_{\min})$
  and 
  \begin{equation*}
   \omega u\xleftarrow{n\to\infty}\omega u_n=BA(\omega u_n)
   \xrightarrow{n\to\infty}BA(\omega u)
  \end{equation*}
  in $\calH^{0,\gamma}_p(\bz)$. The continuity of 
$B:\calH^{0,\gamma}_p(\bz)\to \calH^{\mu,\gamma+\mu-\epsilon}_p(\bz)$
implies that
   $$\calD(A_{\min})\subset\calD(A_{\max})\cap
     \mathop{\mbox{\Large$\cap$}}_{\varepsilon>0}\,
     \calH^{\mu,\gamma+\mu-\varepsilon}_p(\bz)=:V.$$
  Since $A$ has constant coefficients, $\calD(A_{\max})=\calD(A_{\min})
\oplus \calE$,
  where $\calE$ has zero intersection with 
  $\mathop{\mbox{$\cap$}}\limits_{\varepsilon>0}\,
  \calH^{\mu,\gamma+\mu-\varepsilon}_p(\bz)$, see \cite{Lesc}, Proposition
  1.3.11. From this we immediately obtain $\calD(A_{\min})=V$. 
 \end{proof}

%%%%%%%%%%%%%%%%%%%%%%%%%%%%%%%%%%%%%%%%%%%%%%%%%%%%%%%%%%%%%%%%%%%%%

\subsection{The maximal extension}\label{section1.4} 

 Before characterizing the domain of the maximal extension, we shall discuss a
 certain type of operators, namely those of the form 
  \begin{equation}\label{shift}
   G=\omega\left(\op_M^{\gamma_1-\frac{n}{2}}(g)-
   \op_M^{\gamma_2-\frac{n}{2}}(g)\right):
   \cicomp(\partial\bz^\wedge)\longrightarrow\ci(\intb).
  \end{equation}
 Here, $\partial\bz^\wedge:=\rz_+\times\partial\bz$. Moreover, $g$ is a
 meromorphic  Mellin symbol with asymptotic type $P$ as in 
\cite{Seil2} or Section \ref{appendix5}. 
%% also we could take as domain
%% $\calS^\infty(\partial\bz^\wedge)$ or
%% $\cicomp({]a,b[}\times\partial\bz)$ with arbitrary $0<a<b<\infty$ without
%% changing the image of $G$). 
 Let 
  \begin{equation}\label{shift2}
   \smsum_{k=0}^{n_p}R_{pk}(z-p)^{-(k+1)},\qquad 
   R_{pk}\in L^{-\infty}(\partial\bz),
  \end{equation}
 denote the principal part of $g$ around $p\in\pi_\cz P$. Recall that the
 $R_{pk}$ have finite rank by definition. 

 \begin{lemma}\label{shift3}
  Let $G$ be as in {\rm\eqref{shift}} with $\gamma_1<\gamma_2$, and $R_{pk}$
  as in {\rm\eqref{shift2}}. Then $G$ is of finite rank and, for
  $u\in\cicomp(\partial\bz^\wedge)$,   
   $$(Gu)(t,x)=\omega(t)
     \smsum_{\substack{p\,\in\,\pi_\cz P \\ 
                     -\gamma_2<\text{\rm Re}\,p-\frac{n+1}{2}<-\gamma_1}}
     \smsum^{n_p}_{l=0}\zeta_{pl}(u)(x)\,t^{-p}(\log t)^l$$
  with the linear maps $\zeta_{pl}:\cicomp(\partial\bz^\wedge)\to
  \text{\rm im}\,R_{pl}+\ldots+\text{\rm im}\,R_{pn_p}\subset\ci(\partial\bz)$
  given by 
   $$\zeta_{pl}(u)(x)=\smsum_{k=l}^{n_p}\frac{(-1)^{l}}{l! (k-l)!}
     R_{pk}\frac{\partial^{k-l}}{\partial z^{k-l}}(\calM u)(p,x),$$
  where $\calM=\calM_{t\to z}$ denotes the Mellin transform. 
 \end{lemma}
 
 The proof is a straightforward consequence of the residue theorem, since 
  $$(Gu)(t,x)=\Big(\int_{\Gamma_{\frac{n+1}{2}-\gamma_1}}-
    \int_{\Gamma_{\frac{n+1}{2}-\gamma_2}}\Big) 
    t^{-z}g(z)(\calM u)(z,x)\,\dbar z=
    \int_{\calC}t^{-z}g(z)(\calM u)(z,x)\,\dbar z$$
 with a path $\calC$ simply surrounding the poles of $g$ in the strip 
 $\frac{n+1}{2}-\gamma_2<\re z<\frac{n+1}{2}-\gamma_1$. For the detailed
 calculations and an expression for $\mbox{\rm rank}\,G$ see
 \cite{Lesc}. 
 The residue theorem also implies that we could
replace $g$ by $g+h$ for any $h\in
 M^\mu_O(\partial\bz)$ without changing $G$.

 \begin{remark}\label{split}
  Let $\gamma$ with $\gamma_1<\gamma<\gamma_2$ be given and let  
  $G_j=\omega\left(\op_M^{\gamma_j-\frac{n}{2}}(g)-
  \op_M^{\gamma-\frac{n}{2}}(g)\right)$. Then $G=G_1-G_2$ and  
   $$\mbox{\rm im}\,G=\mbox{\rm im}\,G_1\oplus\mbox{\rm im}\,G_2.$$
 \end{remark}

 In fact, by the previous lemma, the images on the right-hand side have
 trivial intersection, and $G_2u_2$ only depends on finitely many Taylor
 coefficients of the Mellin transform $\calM u_2$ in the poles of $g$ lying
 in the strip $\frac{n+1}{2}-\gamma_2<\re z<\frac{n+1}{2}-\gamma$. The
 analogous statement holds for $G_1u_1$. Then the result follows from the
 following observation: Given finitely many points $p_1,\ldots,p_N\in\cz$
 and, in each of these, a finite number of Taylor coefficients, there
 exists a $u\in\cicomp(\rz_+)$ such that the Taylor expansion of $\calM u$
 in each $p_j$ starts with these prescribed values. 

 Now let $A$ be as in \eqref{coneoperator} and set
  \begin{equation}\label{taylor}
   f_l(z)=\frac{1}{l!}\smsum_{j=0}^\mu(d^l_t a_j)(0)z^j,\qquad
   l=0,\ldots,\mu-1. 
  \end{equation}
 In particular, $f_0=\sigma_M^\mu(A)$ is the conormal symbol of $A$. Due to 
 the $\bz$-ellipticity of $A$, $f_0$ is meromorphically invertible and 
 $f^{-1}_0$ can be written as the sum of a meromorphic Mellin symbol and a 
 holomorphic symbol in $M^{-\mu}_O(\partial\bz)$ (see Theorem 6 in Section
 2.3.1 of \cite{Schu2}). We now define recursively 
  \begin{equation}\label{recursiv}
    g_0=f_0^{-1},\qquad
    g_l=-(T^{-l}f^{-1}_0)\smsum_{j=0}^{l-1}(T^{-j}f_{l-j})g_j,\qquad
    l=1,\ldots,\mu-1,
  \end{equation}
 with $T^\sigma$, $\sigma\in\rz$, acting on meromorphic functions by 
 $(T^\sigma f)(z)=f(z+\sigma)$.  

 Moreover, choose an $\varepsilon>0$ so small that every pole $p$ of one of the
 symbols $g_0,\ldots,g_{\mu-1}$, either lies on one of the lines 
 $\Gamma_{\frac{n+1}{2}-\gamma-\mu+k}$, $k=0,\ldots,\mu$, or has a distance to
 each of these lines which is larger then $\varepsilon$. 

 \begin{definition}\label{asymptotic4}
  Let $g_0,\ldots,g_{\mu-1}$ be as in \eqref{recursiv} and $\varepsilon>0$ 
  as described before. Then we set 
   $$\calE=\calE^\gamma_A=\mbox{\rm im}\,G_0+\ldots+\mbox{\rm im}\,G_{\mu-1},$$
  where the operators 
  $G_k=\sum\limits_{l=0}^kG_{kl}:
  \cicomp(\partial\bz^\wedge)\to \calH^{\infty,\gamma+\varepsilon}_p(\bz)$
%\calC^{\gamma,\infty+\varepsilon}(\bz)$ 
  are defined by
   $$G_0=G_{00}=\omega\left(\op_M^{\gamma+\mu+\varepsilon-1-\frac{n}{2}}(g_0)-
     \op_M^{\gamma+\mu-\varepsilon-\frac{n}{2}}(g_0)\right),$$
  and if $1\le k\le\mu-1$, $0\le l\le k$, 
   $$G_{kl}=\omega\, t^l
        \left(\op_M^{\gamma+\mu+\varepsilon-k-1-\frac{n}{2}}(g_l)-
        \op_M^{\gamma+\mu+\varepsilon-k-\frac{n}{2}}(g_l)\right).$$
 \end{definition} 
 
 The space $\calE$ is a finite-dimensional subspace of 
 $\calC^{\infty,\gamma+\varepsilon}(\bz)$ and consists of functions
 of the form 
  \begin{equation}\label{asymptotic5}
   u(t,x)=\omega(t)\,\smsum_{j=0}^N\smsum_{k=0}^{l_j} 
   u_{jk}(x)\,t^{-q_j}\,\log^k t
  \end{equation}
 with smooth functions $u_{jk}\in\ci(\partial\bz)$ and complex numbers $q_j$
 with 
  \begin{equation}\label{streifen}
   \frac{n+1}{2}-\gamma-\mu\le\re q_j<\frac{n+1}{2}-\gamma. 
  \end{equation}
 Note that in case $A$ has constant coefficients we have, due to 
 Remark \ref{split}, 
  \begin{equation}\label{constant}
   \calE=
    \mbox{\rm im}\,G_{00}\oplus\ldots\oplus\mbox{\rm im}\,G_{(\mu-1)0}=
    \mbox{\rm im}\,\omega\left(
     \op_M^{\gamma+\mu-\varepsilon-\frac{n}{2}}(\sigma_M^\mu(A)^{-1})-
     \op_M^{\gamma+\varepsilon-\frac{n}{2}}(\sigma_M^\mu(A)^{-1})\right);
  \end{equation}
 in particular, we have twice strict inequality `$<$' in
 \eqref{streifen}. For $A$ having non-constant coefficients, equality in
 \eqref{streifen} is possible, see Example \ref{abcd}, below. 
 
 \begin{proposition}\label{basis}
  For any $0\le k\le \mu-1$ let $u_{k1},\ldots,u_{kn_k}\in
  \cicomp(\partial\bz^\wedge)$ be chosen such that 
  $\{G_{k0}u_{kj}\st 1\le j\le n_k\}$ is a basis of $\mbox{\rm im}\,G_{k0}$. 
  Then 
   $$\{G_ku_{kj}\st 0\le k\le\mu-1,\, 1\le j\le n_k\}\subset\calE$$
  is a set of linearly independent functions such that 
   $$\mbox{\rm span}\{G_ku_{kj}\st 0\le k\le\mu-1,\, 1\le j\le n_k\}\cap
     \calD(A_{\min})=\{0\}.$$
  In particular, 
   $$\mbox{\rm dim}\,\calE\ge\mbox{\rm dim\,im}\,\omega
     \left(\op_M^{\gamma+\varepsilon-\frac{n}{2}}(\sigma^\mu_M(A)^{-1})-
     \op_M^{\gamma+\mu-\varepsilon-\frac{n}{2}}(\sigma^\mu_M(A)^{-1})\right).$$
  We have equality at least in the cases where $A$ has constant coefficients
  near the boundary or $\sigma^\mu_M(A)^{-1}$ has no pole on the line 
  $\re z=\frac{n+1}{2}-\gamma-\mu$. 
 \end{proposition}
 
 \begin{proof}
  Let $\alpha_{jk}\in\cz$ and 
   $$\smsum_{k=0}^{\mu-1}\smsum_{j=1}^{n_k}\alpha_{jk}G_ku_{kj}=u\in
     \calD(A_{\min})\subset\calH^{0,\gamma+\mu-\varepsilon}_p(\bz).$$
  Setting $l=\mu-1$, we obtain 
   $$\smsum_{j=1}^{n_l}\alpha_{lj}G_{l0}u_{lj}=u-
     \smsum_{k=0}^{l-1}\smsum_{j=1}^{n_k}\alpha_{jk}G_ku_{kj}-
     \smsum_{j=1}^{n_l}\alpha_{lk}(G_l-G_{l0})u_{lj}.$$
  The right-hand side belongs to $\calH^{0,\gamma+\mu-\varepsilon}_p(\bz)+
  \calH^{0,\gamma+1+\varepsilon}_p(\bz)$. The intersection of this space 
  with $\mbox{\rm im}\,G_{l0}$ is trivial, hence 
  $\sum\limits_{j=1}^{n_l}\alpha_{lj}G_{l0}u_{lj}=0$. Therefore
  $\alpha_{lj}=0$ for all $1\le j\le n_l$, since the $G_{l0}u_{lj}$ are
  linearly independent by assumption. Iterating this process (i.e. taking
  $l=\mu-2$, $l=\mu-3$, etc.), we see that all $\alpha_{jk}$ must equal zero. 
  
  If $A$ has constant coefficients, the result on the dimension follows from
  \eqref{constant}; the second identity in \eqref{constant} is always true and
  yields the lower bound for $\mbox{\rm dim}\,\calE$. The remaining claim we
  shall obtain as a by-product of the following theorem. 
 \end{proof}
  
 \begin{theorem}\label{maxdomain}
  The domain of the maximal extension $A_{\max}=A^{\gamma,p}_{\max}$ of $A$
  from \eqref{unbounded} is 
   $$\calD(A_{\max})=\calD(A_{\min})+\calE$$
  with $\calE$ from Definition {\rm\ref{asymptotic4}}. Recall that $\calE$ does
  not depend on $1<p<\infty$. The sum is direct at
  least in the cases where $A$ has constant coefficients near the boundary or
  $\sigma^\mu_M(A)^{-1}$ has no pole on the line 
  $\re z=\frac{n+1}{2}-\gamma-\mu$. In any case, 
   $$\calD(A_{\min})\cap\calE\subset
     \mbox{\rm im}\,\omega
     \left(\op_M^{\gamma+\mu-\varepsilon-\frac{n}{2}}(\sigma^\mu_M(A)^{-1})-
     \op_M^{\gamma+\mu+\varepsilon-\frac{n}{2}}(\sigma^\mu_M(A)^{-1})\right).$$
  Consequently, any closed extension $\underline{A}$ in
  $\calH^{0,\gamma}_p(\bz)$ is given by the action of $A$ on a domain 
   \begin{equation}\label{domain}
    \calD(\underline{A})=\calD(A_{\min})+\underline{\calE},\qquad 
    \underline{\calE}\text{ subspace of }\calE.
   \end{equation} 
 \end{theorem}
 
 \begin{proof}
  For a cut-off function $\widetilde{\omega}$ supported
  sufficiently close to zero, the operator $\widetilde{A}=\widetilde{\omega}\,
  t^{-\mu}$\, $\sum\limits_{j=0}^\mu a_j(0)(-t\partial_t)^j+
  (1-\widetilde{\omega})A$ is still $\bz$-elliptic, and 
  $\calD(\widetilde{A}_{\min})=\calD(A_{\min})$ by \cite{Lesc}, Proposition
  1.3.12. Moreover, 
   $$\calD(\widetilde{A}_{\max})=\calD(\widetilde{A}_{\min})\oplus
     \widetilde{\calE}$$
  with $\widetilde{\calE}$ given by the right-hand side of \eqref{constant},
  and 
   $$\mbox{\rm dim}\,\calD(A_{\max})\big/\calD(A_{\min})=
     \mbox{\rm dim}\,\calD(\widetilde{A}_{\max})\big/
     \calD(\widetilde{A}_{\min})=\mbox{\rm dim}\,\widetilde{\calE}.$$
  The latter statements are due to \cite{Lesc}, Proposition 1.3.11,
  Corollary 1.3.17. By Proposition
  \ref{basis} it therefore suffices to prove that $\calE\subset\calD(A_{\max})$.
  In fact, we shall show that $\mbox{\rm im}\,G_k$ belongs to $\calD(A_{\max})$
  for any $k$. Since this is easy to see for $k=0$, we shall only consider the
  case $k\ge1$. With the holomorphic Mellin symbols $f_n$ from
  \eqref{taylor} write
   $$A=\widetilde{\omega}\,t^{-\mu}\smsum_{j=0}^{\mu-1}t^j\op_M(f_j)+
   t^\mu A^\prime$$
  for a $\mu$-th order cone differential operator $A^\prime$. Taking into
  account that $\mbox{\rm im}\,G_{kl}$ is a subset of 
  %$\calC^{\infty,\gamma+\mu+\varepsilon-k+l-1}(\bz)$ 
  $\calH^{\infty,\gamma+\mu+\varepsilon-k+l-1}_p(\bz)$ we thus obtain for 
  $u\in\cicomp(\partial\bz^\wedge)$ 
   $$A(G_ku)\in\calH^{0,\gamma}_p(\bz)\iff
     \widetilde{\omega}\smsum_{j=0}^k\smsum_{l=0}^{k-j}
     t^j\op_M(f_j)(G_{kl}u)\in\calH^{0,\gamma+\mu}_p(\bz).$$
  Choosing $\widetilde{\omega}$ with
  $\widetilde{\omega}\omega=\widetilde{\omega}$, using the elementary rule 
  \begin{equation}\label{kommu}
\op_M(f)t^\sigma\op_M^\delta(g)=t^\sigma\op_M^\delta((T^{-\sigma}f)g),
\end{equation}
  and rearranging the order of summation, we see that 
  $A(G_ku)\in\calH^{0,\gamma}_p(\bz)$ if and only if   
   $$\widetilde{\omega}\smsum_{j=0}^kt^j\smsum_{l=0}^j\left(
     \op_M^{\gamma+\mu+\varepsilon-k-1-\frac{n}{2}}((T^{-l}f_{j-l})g_l)-
     \op_M^{\gamma+\mu+\varepsilon-k-\frac{n}{2}}((T^{-l}f_{j-l})g_l)\right)(u)
     \in\calH^{0,\gamma+\mu}_p(\bz).$$
  However, this expression actually equals zero, since by definition of the
  symbols $g_l$, cf.\ \eqref{recursiv}, we have 
   $$\smsum_{l=0}^j(T^{-l}f_{j-l})g_l=\delta_{0j},\qquad 0\le j\le k,$$
  with $\delta_{0j}$ denoting the Kronecker symbol. This shows the claim. 
  
  Let us turn to the remaining claims of the theorem. If $A$ has constant
  coefficients near the boundary the intersection of $\calD(A_{\min})$ and
  $\calE$ is zero by \eqref{constant}. Using the description of
  elements $u$ from $\calE$ given in \eqref{asymptotic5}, and the fact that
  $\calD(A_{\min})\subset\calH^{0,\gamma+\mu-\delta}_p(\bz)$ for any positive
  $\delta$, we see that if $u\in\calD(A_{\min})\cap\calE$ then $u$ is of the
  form \eqref{asymptotic5} with $\re q_j=\frac{n+1}{2}-\gamma-\mu$. 
  Let $v:=t^{-\mu}\,\op_M(\sigma_M^\mu(A))(u)$.
  Inserting the form \eqref{asymptotic5},
  the fact that $v\in\calH^{0,\gamma}_p(\bz)$, implies that
  $v\in\cicomp(\intb)$, in particular, $\widetilde{\omega}v=0$ for a suitable
  cut-off function $\widetilde{\omega}$. We now set
  $B_{\pm\varepsilon}=
  \op_M^{\gamma+\mu\pm\varepsilon-\frac{n}{2}}(\sigma^\mu_M(A)^{-1})\,t^\mu$
  for small $\varepsilon>0$. Then 
   $$0=B_{-\varepsilon}(\widetilde{\omega}v)=
     B_{-\varepsilon}v+B_{-\varepsilon}((1-\widetilde{\omega})v)=
     u+(B_{-\varepsilon}-B_{\varepsilon})((1-\widetilde{\omega})v)+
     B_{\varepsilon}((1-\widetilde{\omega})v).$$
  We conclude from Lemma \ref{shift3} that the last term is zero and 
  $u\in\mbox{\rm im}(B_{-\varepsilon}-B_{\varepsilon})$. 
  In particular, the intersection $\calD(\amin)\cap \calE$ is trivial if $\sigma_M^\mu(A)^{-1}$ has no
  pole on the line $\re z=\frac{n+1}{2}-\gamma-\mu$. This then also proves the
  last claim of Proposition \ref{basis} as we announced in the previous proof.
  
  Since $\calE$ is finite-dimensional all the operators $\underline{A}$ in \eqref{domain} are in
  fact closed extensions of $A$. 
 \end{proof}
 
 From Theorem \ref{maxdomain} and Proposition \ref{mindomain} one might
 conjecture that 
  $$\calD(A_{\min})=\calH^{\mu,\gamma+\mu}_p(\bz)\oplus
    \mbox{\rm im}\,\omega\left(
    \op_M^{\gamma+\mu-\varepsilon-\frac{n}{2}}(\sigma^\mu_M(A)^{-1})- 
    \op_M^{\gamma+\mu+\varepsilon-\frac{n}{2}}(\sigma^\mu_M(A)^{-1})\right)$$
 in case $\sigma^\mu_M(A)^{-1}$ has a pole with real part
 $\frac{n+1}{2}-\gamma-\mu$. This however is {\em not} true, for then
 $A:\calH^{\mu,\gamma+\mu}_p(\bz)\to\calH^{0,\gamma}_p(\bz)$ would be a
 Fredholm operator, thus have invertible conormal symbol (by equivalence of
 ellipticity and Fredholm property in the cone algebra, cf.\ \cite{ScSe}, 
 Theorem 3.13). 

 \begin{example}\label{abcd}
  Let $\bz$ have dimension {\rm2} and boundary $\partial\bz=S^1$ {\rm(}the unit 
  sphere{\rm)}.  Define%With a cut-off function $\omega$ let 
   $$A=\omega\,t^{-2}\left\{\mbox{$\frac{1}{4}$}
     ((t\partial_t)^2-t(t\partial_t))+\Delta_\partial\right\}+
     (1-\omega)\Delta,$$
  where $\Delta_\partial$ is the standard Laplacian on $S^1$ and $\Delta$ some 
  Laplacian on $\intb$.
% Then 
%   $$\sigma^2_\psi(A)(y,\eta)=
%     \omega(t)t^{-2}(\mbox{$\frac{1}{4}$}(t\tau)^2+|\xi|^2)+
%     (1-\omega)(y)|\eta|^2,$$
%  where, near the boundary, $y=(t,x)$ with covariable $\eta=(\tau,\xi)$. Hence
  Clearly, $A$ is elliptic; it is $\bz$-elliptic, since the rescaled symbol 
is $-\tau^2/4-|\xi|^2$. 
The conormal symbol  
   $$\sigma^2_M(A)(z)=\mbox{$\frac{1}{4}$}z^2+\Delta_\partial$$
  has the non-bijectivity points $z\in2\gz$, since the spectrum of
  $\Delta_\partial$ is $\{-k^2\st k\in\nz_0\}$. 

Considering $A$ as an unbounded
  operator in $\calH^{0,0}_p(\bz)$
%   $$A:\cicomp(\intb)\subset\calH^{0,0}_p(\bz)\longrightarrow
%     \calH^{0,0}_p(\bz)$$
  we obtain from Proposition {\rm\ref{mindomain}} that
  $\calD(A_{\min})=\calH^{2,2}_p(\bz)$. From Theorem {\rm\ref{maxdomain}} and
  Proposition {\rm\ref{basis}} we conclude that 
  $\calD(A_{\max})=\calH^{2,2}_p(\bz)\oplus\calE$ with a two dimensional space
  $\calE$. Actually, by direct computation,  
   $$\calE=\mbox{\rm span}\{\omega,\,\omega(t+\log t)\}.$$
  In the  notation  of Definition {\rm\ref{asymptotic4}}, the functions
  $\omega$,  $\omega\log t$ are generated by $G_{10}$, while $\omega t$ comes
  from $G_{11}$. The operator $G_{00}$ equals zero, since $\sigma_M^2(A)^{-1}$
  has no poles in the strip $-1\le\re z<0$. 
 \end{example}

%%%%%%%%%%%%%%%%%%%%%%%%%%%%%%%%%%%%%%%%%%%%%%%%%%%%%%%%%%%%%%%%%%%%%%

\subsection{The model cone operator}\label{section1.5}
 
 Freezing the coefficients of $A$ at the boundary leads to the differential
 operator 
  \begin{equation}\label{modelcone}
   \widehat{A}=t^{-\mu}\smsum_{j=0}^\mu a_j(0)(-t\partial_t)^j
  \end{equation}
 on the infinite half-cylinder $\partial\bz^\wedge=\rz_+\times\partial\bz$. We
 shall refer to this operator as the {\em model cone operator} of $A$. Let us
 first introduce a suitable scale of Sobolev spaces the operator $\widehat{A}$
 acts in. 
 
 To this end let $\partial \bz=X_1\cup \ldots\cup X_J$ be an open covering of 
 $\partial\bz$; let $\kappa_j:X_j\to U_j$ be coordinate maps and
 $\{\varphi_1,\ldots , \varphi_J\}$ a subordinate partition of unity.

 Given a function $u=u(t,x) $ on $\rz\times\partial\bz$, we shall say that
 $u\in H_{p,{\rm cone}}^s(\rz\times\partial\bz)$ provided that, for each $j$,
 the function
  $$v(t,y) = \varphi_j(x)u(t,x),\quad x= \kappa^{-1}_j(y/[t]),$$
 is an element of $H_p^s(\rz\times \rz^n)$ $($we consider the right-hand side
 to be zero for $x\notin U_j)$. 
 In other words: $\varphi_j u$  is the pull-back of a function in
 $H_p^s(\rz^{n+1})$ under the composition of the maps
  $${\rm id}\times \kappa_j:\rz\times
    X_j\ni (t,x)\mapsto (t,\kappa_j(x))\in \rz\times U_j$$
 and
  $$\chi:\rz\times U_j\ni(t,y){\mapsto}(t,[t]y)\in\rz^{n+1},$$ 
 so that the definition extends to distributions in the usual way for
 $s\in\rz$, $1<p<\infty$.

 \begin{definition}
  ${\mathcal K}^{s,\gamma}_p(\partial\bz^\wedge)$ is the space of all
  distributions $u\in H^s_{p,{\rm loc}}(\rz_+\times \partial\bz)$ such that,
  for an arbitrary cut-off function $\omega$,
   $$\omega u \in {\mathcal H}^{s,\gamma}_p(\bz) \quad{\rm and }
     \quad (1-\omega) u\in H^s_{p,{\rm cone}}(\rz\times \partial\bz).$$
 \end{definition}
 
 $\widehat{A}$ acts continuously in this scale of Sobolev spaces, 
  $$\widehat{A}:\calK^{s,\gamma}_p(\partial\bz^\wedge)\longrightarrow 
    \calK^{s-\mu,\gamma-\mu}_p(\partial\bz^\wedge),\qquad
    s,\gamma\in\rz,\;1<p<\infty.$$
 We shall now consider the model cone operator as an unbounded operator,
 namely 
  $$\widehat{A}:\cicomp(\partial\bz^\wedge)\subset
    \calK^{0,\gamma}_p(\partial\bz^\wedge)\longrightarrow 
    \calK^{0,\gamma}_p(\partial\bz^\wedge).$$
 %If $A$ satisfies condition {\rm(1)} of Section \ref{section2.1} below, 
We first show that the
 domains of the closed extensions of $\widehat{A}$ can be read off from those of $A$
provided $A$ satisfies a mild additional ellipticity condition. In
 analogy to Theorems \ref{mindomain} and \ref{maxdomain} we have: 

 \begin{proposition}\label{elmar}
  If $A$ satisfies condition {\rm(E1)} of Section {\rm\ref{section2.1}}, then 
   $$\calD(\widehat{A}_{\min})=
     \Big\{u\in\mathop{\mbox{\Large$\cap$}}_{\varepsilon>0}\,
     \calK^{\mu,\gamma+\mu-\varepsilon}_p(\partial\bz^\wedge)\st 
     \widehat{A}u\in\calK^{0,\gamma}_p(\bz)\Big\}.$$
  This simplifies to 
  $\calD(\widehat A_{\rm min})=\calK^{\mu,\gamma+\mu}_p(\partial \bz^\wedge)$ 
  if and only if $\sigma_M^\mu(A)$ is invertible on the line
  $\Gamma_{\frac{n+1}2-\gamma-\mu}$.
  
  If $\calE$ is the space defined in \eqref{constant}, then 
   $$\calD(\widehat{A}_{\max})=\calD(\widehat{A}_{\min})\oplus\calE.$$     
  Hence, any closed extension $\underline{\widehat{A}}$ of $\widehat{A}$ is
  given by the action of $\widehat{A}$ on a domain 
   $$\calD(\underline{\widehat{A}})=\calD(\widehat{A}_{\min})\oplus
     \underline{\calE},\qquad\underline{\calE}\text{ subspace of }
     \calE.$$
 \end{proposition}

 \begin{proof}
  Let $\widetilde P = \sum _{j=0} ^\mu a_j(0)(-\partial_t)^j-\eta^\mu$. This
  is a non-degenerate parameter-dependent differential operator with
  coefficients independent of $t$. It follows from {\rm(E1)} in Section
  \ref{section2.1} that the parameter-dependent principal symbol 
   $$\widetilde p_0(x,\xi,\tau,\eta)=\smsum_{j=0}^\mu
     \sigma_\psi^{\mu-j}(a_j)(0,x,\xi)(-i\tau)^j-\eta^\mu$$
  of $\widetilde{P}$ is parameter-elliptic. Hence there exist  parameter-dependent
  symbols $\widetilde q_0$ of order $-\mu$ and $r_1,r_2$ of order $-1$ such that 
   \begin{equation}\label{paramell}
    \widetilde p_0 \widetilde q_0 = 1+r_1, \quad  \widetilde q_0 \widetilde p_0 = 1+r_2.
   \end{equation}
  The operator $P= t^\mu(\widehat A-\eta^\mu)$ has the principal  symbol
   $$\sigma_\psi^\mu(P)(x,t,\xi,\tau,\eta) = \widetilde p_0(x,\xi,t\tau,t\eta).$$
  Under the push-forward induced by $T:=\chi\circ({\rm id}\times\kappa_j)$ the
  operator $P$ transforms into a weighted $SG$-pseudodifferential operator 
  of order $(\mu,\mu)$; modulo terms of order $(\mu-1,\mu-1)$ its
  $SG$-symbol is given by the push-forward of $\sigma_\psi^\mu(P)$. Indeed, for
  a differential operator this is a simple calculation, a proof for the general
  pseudodifferential case is given in \cite{ScSc}, Theorem 3.8; details on
  $SG$-symbols can be found in \cite{Schr}. Now equation \eqref{paramell}
  implies that the push-forward of $\sigma_\psi^\mu(P)$ is
  $SG$-parameter-elliptic if we restrict to a subset of 
  $\rz_+\times\partial\bz$ away from the boundary, say to $\{t\ge1\}$.
  Hence, on this set, there is a parameter-dependent $SG$-parametrix $S$ of
  order $(-\mu,-\mu)$ to the push-forward $T_*P$ of $P$ (i.e.\ we have 
  $S\circ T_*P =I+R$, where $R$ is an integral operator with a rapidly
  decreasing kernel). As the operator of multiplication
  by $t^\mu$ remains unchanged under the
  push-forward and is an $SG$-operator of order $(0,\mu)$, $S\circ t^\mu$ is an
  $SG$-parametrix of order $(-\mu,0)$ to the push-forward of $A-\eta^\mu$. 

 We now can describe the domain away from the boundary: Given  
  $u\in \calK^{0,\gamma}_p$, a cut-off function $\omega$
  equal to $1$ in $\{t\le 1\}$, and a function $\varphi_j$ in the partition of
  unity on $\partial\bz$,
   \begin{equation}\label{sgparametrix}
    S\circ t^\mu\circ
    T_*\left((\widehat A-\eta^\mu)(1-\omega)\varphi_ju \right)
    =T_*(1-\omega)\varphi_ju +RT_*(1-\omega)\varphi_ju
   \end{equation}
  Let $u\in \calD(\widehat A_{\rm max})$, so that 
  $\widehat Au \in \calK^{0,\gamma}_p$. Standard elliptic regularity
  implies that $u\in H^\mu_{p,{\rm loc}}$. Hence $\widehat
  A(1-\omega)\varphi_ju\in \calK^{0,\gamma}_p$: Outside a compact set, it 
  coincides with the function
  $(1-\omega)\varphi_j\widehat Au$ whose  push-forward
  via $T$  belongs to $H_p^0(\rz^{n+1})$. In view of the fact that  $St^\mu:
  H_p^0(\rz^{n+1}) \to H_p^\mu(\rz^{n+1})$ is bounded, we deduce from
  \eqref{sgparametrix} that $(1-\omega)u\in 
  H^\mu_{p,{\rm cone}}(\rz\times\partial\bz)$. On the other hand, we trivially
  have $u\in \calD(\widehat A_{\rm min})$ for every $u$ in 
  $H^\mu_{p,{\rm cone}}(\rz\times \partial\bz)$ supported in $\{t\ge1\}$.
  As a consequence, the domains of all extensions of $\widehat A$ coincide
  with $H^\mu_{p,{\rm cone}}(\rz\times \partial\bz)$ away from $\{t=0\}$.
  
  Close to $\{t=0\}$, the analysis is the same as in the standard
  case. This completes the argument. 
 \end{proof}
%%%%%%%%%%%%%%%%%%%%%%%%%%%%%%%%%%%%%%%%%%%%%%%%%%%%%%%%%%%%%%%%%%%%%%
%%%%%%%%%%%%%%%%%%%%%%%%%%%%%%%%%%%%%%%%%%%%%%%%%%%%%%%%%%%%%%%%%%%%%%
 
\section{Structure of the resolvent}\label{section2}
 
 Let us now come to the main objective of this paper. We shall consider a
 closed extension of a cone differential operator and give conditions that
 ensure that its resolvent exists in a given sector $\Lambda$
 (up to finitely many exceptional points). We shall describe the structure of
 this resolvent in terms of a class of parameter-dependent cone
 pseudodifferential operators. 
 
 Before considering the resolvent, we want to investigate the inverse of a
 given closed extension. This is a simpler problem but already illustrates some
 of the structures we shall see in the discussion of the resolvent.  
 We again refer to \cite{Seil2} and the appendix for basic notions of the 
 cone calculus.  
%%%%%%%%%%%%%%%%%%%%%%%%%%%%%%%%%%%%%%%%%%%%%%%%%%%%%%%%%%%%%%%%%%%%%%
 
\subsection{The inverse of a closed extension}\label{section1.3}
   
 Let $A$ be a cone differential operator and assume that  
  $$\underline{A}:\calD(\underline{A})=\calH^{\mu,\gamma+\mu}_p(\bz)\oplus 
    \underline{\calE}\longrightarrow\calH^{0,\gamma}_p(\bz)$$
 is bijective for a fixed choice of $\gamma$ and $p$. 
 Since $\underline{\calE}$ is finite-dimensional,  
 $A:\calH^{\mu,\gamma+\mu}_p(\bz)\to\calH^{0,\gamma}_p(\bz)$ is a Fredholm 
 operator. According to  \cite{ScSe}, Theorem 3.13, this is 
%we have shown that the Fredholm property is 
 equivalent to the ellipticity of  $A$,  i.e.\ $A$ is $\bz$-elliptic and 
 the conormal symbol is invertible on the line 
 $\Gamma_{\frac{n+1}{2}-\gamma-\mu}$. In other words, $A$ is an 
 elliptic element of the cone algebra $C^\mu(\gamma+\mu,\gamma,k)$ for 
 any $k\in\nz$.
% The cone algebra on $\bz$ was introduced by Schulze; for its 
% definition we refer to \cite{Schu1} (and, concerning notation we use here, to
% the corresponding definitions of the parameter-dependent version given in the
% appendix). 
 
 Let us now set $\underline{\calF}=A(\underline{\calE})$. This space is 
 finite-dimensional and 
 $\calH^{0,\gamma}_p(\bz)=
 A(\calH^{\mu,\gamma+\mu}_p(\bz))\oplus\underline{\calF}$.

 \begin{lemma}\label{lemmaa}
  There exists an asymptotic type $Q\in\mbox{\rm As}(\gamma,k)$, $k\in\nz$ 
  arbitrarily large, such that 
  $\underline{\calF}=A(\underline{\calE})\subset
  \calC^{\infty,\gamma}_Q(\bz).$ 
 \end{lemma}

 \begin{proof}
  Let $u$ be of the form \eqref{asymptotic5}, and choose a cut-off function
  $\widetilde{\omega}$ with $\widetilde{\omega}\omega=\widetilde{\omega}$. Then, with $A$
  as in \eqref{coneoperator}, 
   $$(Au)(t)=\widetilde{\omega}t^{-\mu}\smsum_{j=0}^\mu\smsum_{q,k}
     a_j(t)(t\partial_t)^j(c_{qk}t^{-q}\log^k t)+
     (1-\widetilde{\omega})(t)(Au)(t).$$
  The second term belongs to $\cicomp(\intb)$. Now a Taylor expansion of the
  coefficients $a_j$ in $t$ at 0 shows the claim.
 \end{proof}
 
 Let us now denote by 
  $$\pi_{\underline{\calE}}:\calH^{\mu,\gamma+\mu}_p(\bz)\oplus
    \underline{\calE}\longrightarrow\underline{\calE},\qquad 
    \pi_{\underline{\calF}}:A(\calH^{\mu,\gamma+\mu}_p(\bz))\oplus
    \underline{\calF}\longrightarrow\underline{\calF}$$ 
 the canonical projections, and let 
  $$B=(1-\pi_{\underline{\calE}})\underline{A}^{-1}:
    \calH^{0,\gamma}_p(\bz)\longrightarrow\calH^{\mu,\gamma+\mu}_p(\bz)$$
 be a left-inverse of 
 $A:\calH^{\mu,\gamma+\mu}_p(\bz)\to\calH^{0,\gamma}_p(\bz)$. 

 \begin{lemma}\label{lemmab}
  $1-AB=\pi_{\underline{\calF}}$ belongs to 
  $C_G(\gamma,\gamma,k)$ for arbitrarily large $k\in\nz$. 
 \end{lemma}
 \begin{proof}
  By construction, it is clear that $1-AB=\pi_{\underline{\calF}}$. Let 
  $w_1,\ldots,w_N$ be a basis of $\underline{\calF}$. Then we can write 
  $\pi_{\underline{\calF}}(\cdot)=\sum_{j=1}^Nl_j(\cdot)w_j$
  with functionals $l_j$ on $\calH^{0,\gamma}_p(\bz)$. By duality 
  there exist $v_j\in\calH^{0,-\gamma}_{p'}(\bz)$ such that 
  $l_j(\cdot)=\skp{\cdot}{v_j}_{\calH^{0,0}_2(\bz)}$. Then, for all 
  $u\in\calH^{\mu,\gamma+\mu}_p(\bz)$, 
   $$0=\pi_{\underline{\calF}}(Au)=
     \smsum_{j=1}^N\skp{Au}{v_j}_{\calH^{0,0}_2(\bz)}w_j=
     \smsum_{j=1}^N\skp{u}{A^*v_j}_{\calH^{0,0}_2(\bz)}w_j,$$
  where $A^*$ denotes the formal adjoint of $A$, which belongs 
  to $C^\mu(-\gamma,-\gamma-\mu,k)$. Hence $v_j\in\text{ker}\,A^*$ for 
  $j=1,\ldots,N$. Since with $A$ also $A^*$ is an elliptic cone operator, 
  $\text{ker}\,A^*\subset \calC^{\infty,-\gamma}_{Q^\prime}(\bz)$ for some 
  asymptotic type $Q^\prime\in\mbox{\rm As}(-\gamma,k)$ by elliptic 
  regularity, cf.\ \cite{Schu2}, Theorem 8 in Section 2.2.1. Thus
  $\pi_{\underline{\calF}}$ has a kernel in
  $\underline{\calF}\otimes\calC^{\infty,-\gamma}_{Q^\prime}(\bz)$, and 
  therefore is a Green operator. 
 \end{proof}

 \begin{proposition}\label{propa}
  $B$ is an element of $C^{-\mu}(\gamma,\gamma+\mu,k)$ for arbitrarily 
  large $k\in\nz$. 
 \end{proposition}

 \begin{proof}
  Since $A$ is an elliptic element of $C^\mu(\gamma+\mu,\gamma,k)$ as 
  we have seen above, there exists a parametrix 
  $C\in C^{-\mu}(\gamma,\gamma+\mu,k)$, i.e. 
  $$AC-1=G_R\in C_G(\gamma,\gamma,k),\qquad 
    CA-1=G_L\in C_G(\gamma+\mu,\gamma+\mu,k).$$ 
  Multiplying these identities from the left, respectively from the right 
  with $B$ yields $B=C-BG_R$ and $B=CAB-G_LB=C-C(1-AB)-G_LB$. Inserting the 
  first equation into the second gives 
   $$B=C-C\pi_{\underline{\calF}}-G_LC+G_LBG_R.$$
  The third term on the right-hand side is a Green operator, since these 
  form an ideal in the cone algebra. The same is true for the 
  second term in view of Lemma \ref{lemmab} and for the forth due to the continuity 
  of $B$ and the mapping properties of Green operators. 
 \end{proof}

 \begin{theorem}\label{propb} 
  $\underline{A}^{-1}=B+G$ for the above 
  $B\in C^{-\mu}(\gamma,\gamma+\mu,k)$ and a suitable  
  $G\in C_G(\gamma,\gamma,k)$ with arbitrarily large $k\in\nz$. 
 \end{theorem}

 \begin{proof}
  We decompose  
  $\underline{A}^{-1}=(1-\pi_{\underline{\calE}})\underline{A}^{-1}+
  \underline{A}^{-1}\pi_{\underline{\calF}}=B+G$. 
  By Proposition \ref{propa}, $B$ is as claimed. From the proof of Lemma 
  \ref{lemmab} we know that $\pi_{\underline{\calF}}$ has an integral 
  kernel in 
  $\underline{\calF}\otimes\calC^{\infty,-\gamma}_{Q^\prime}(\bz)$ for 
  some type $Q^\prime\in\mbox{\rm As}(-\gamma,k)$. Therefore, 
  $\underline{A}^{-1}\pi_{\underline{\calF}}$ has a kernel in 
  $\underline{\calE}\otimes\calC^{\infty,-\gamma}_{Q^\prime}(\bz)$, hence 
  is a Green operator. 
 \end{proof}
 
 As we shall explain in Section \ref{section3.1}, the invertibility is
 independent of $1<p<\infty$. 
 
%%%%%%%%%%%%%%%%%%%%%%%%%%%%%%%%%%%%%%%%%%%%%%%%%%%%%%%%%%%%%%%%%%%%%%

\subsection{Ellipticity assumptions and resolvent analysis}
\label{section2.1}

 Let $A$ be a cone differential operator of order $\mu$ whose coefficients on
 $[0,1]\times\partial\bz$ are independent of $t$ and let 
  $$\underline{A}:\calD(\underline{A})=\calD(\amin)\oplus\underline{\calE}
    \subset\calH^{0,\gamma}_p(\bz)\longrightarrow\calH^{0,\gamma}_p(\bz)$$
 be a closed extension. Moreover, we assume 
 \begin{itemize}
  \item[(E1)]
Both $\sigma_\psi^\mu(A)$ and $\widetilde{\sigma}_\psi^\mu(A)$ have no
   spectrum in the sector $\Lambda$, 
  \item[(E2)] With the above choice of $\underline{\calE}$, the domain 
   $\calD(\underline{\widehat{A}})=
   \calD(\widehat{A}_{\min})\oplus\underline{\calE}$ of the model cone
   operator is invariant under dilations,
  \item[(E3)] The sector $\Lambda\setminus\{0\}$ contains no spectrum of the
   model cone operator
    $$\underline{\widehat{A}}:\calD(\underline{\widehat{A}})\subset
      \calK^{0,\gamma}_p(\partial\bz^\wedge)\longrightarrow
      \calK^{0,\gamma}_p(\partial\bz^\wedge).$$
 \end{itemize}
 
 In \mbox{\rm(E1)} and \mbox{\rm(E3)}, $\Lambda=\Lambda_\theta$ is a closed
 sector in the complex plane containing zero, i.e.
  $$\Lambda_\theta=\{z\in\cz\st |\text{\rm arg}\,z|\ge\theta\}\cup\{0\},$$
 where $0\le\theta<\pi$ and $-\pi\le\text{\rm arg}\,z<\pi$. 

In order to simplify the exposition, we fix a sector $\Sigma$ such that the
mapping  $\eta\mapsto\eta^\mu$ induces a bijection $\Sigma\to\Lambda$. 
Instead of considering $\lambda-A$ for $\lambda\in\Lambda$ we then study
$\eta^\mu-A$ for $\eta\in\Sigma$.

 Condition \mbox{\rm(E2)} means the following: Whenever $u=u(t,x)$ belongs to
 $\calD(\underline{\widehat{A}})$, the same is true for the
 functions $u_\varrho(t,x)=u(\varrho t,x)$, $\varrho>0$. It is easy to see that
 the domain $\calD(\underline{\widehat{A}})$ is invariant under dilations if
 and only if this is true for
 $\cicomp(\partial\bz^\wedge)\oplus\underline{\calE}$. Gil and Mendoza
 \cite{GiMe} call such a domain {\em saturated}. 
 Condition \mbox{\rm(E2)} is always satisfied for the minimal
 extension, the maximal extension, and for extensions with domain equal to 
 $\calD(A_{\max})\cap\calH^{0,\sigma}_p(\bz)$ and $\gamma-\mu<\sigma\le\gamma$.
 For concrete examples see Section \ref{section4} below. 

 \begin{theorem}\label{maintheorem}
  Let $\underline{A}$ satisfy conditions \mbox{\rm (E1), (E2), and
  (E3)} with respect to the sector $\Lambda$. Then 
  \begin{itemize}
   \item[a)] $\underline{A}$ has at most finitely many spectral points in
    $\Lambda$. 
   \item[b)] There exists a parameter-dependent cone pseudodifferential
   operator    
     $$c(\eta)\in C^{-\mu}_O(\Sigma)+
       C^{-\mu}_G(\Sigma;\gamma,\gamma,\varepsilon)$$
    with a certain $\varepsilon>0$, such that 
    $(\eta^\mu-\underline{A})^{-1}=c(\eta)$ for $\eta\in\Sigma$ with $|\eta|$
    sufficiently large. 
  \end{itemize}
 \end{theorem}
 
 For the notation used in part b) of this theorem we refer to the appendix 
 (see Definition \ref{flatcone} and Definition \ref{green}). Part a)
 of the theorem follows from b): Since the domain of $\underline{A}$ is 
 compactly embedded in $\calH^{0,\gamma}_p(\bz)$, $\underline{A}$ has a compact
 resolvent, hence a discrete spectrum. 

%%%%%%%%%%%%%%%%%%%%%%%%%%%%%%%%%%%%%%%%%%%%%%%%%%%%%%%%%%%%%%%%%%%%%%

The proof of Theorem \ref{maintheorem}
is the most technical part of the paper, and it relies on the material in
the appendix. The basic idea is as follows: We write
$a(\eta)=\eta^\mu-A$. This is a parameter-dependent family of cone differential
operators of order $\mu$, i.e.,
  $$a(\eta)\in C^\mu_O(\Sigma)\subset 
    C^\mu(\Sigma;\widetilde{\gamma},\widetilde{\gamma}-\mu,k)$$
 for any $\widetilde{\gamma}\in\rz$ and $k\in\nz$. Assumption
 \mbox{\rm (E1)} ensures that $a(\eta)$ is $\bz$-elliptic with respect
 to the sector $\Sigma$. 
Hence the conormal symbol of $a(\eta)$ (which is that of $-A$) is
meromorphically invertible as shown in \cite{Schu2}, 2.2.1,
Theorem 14 and 2.3.1, Theorem 16.  As we do not
 require $A$ to be conormally elliptic with respect to the weight $\gamma+\mu$,
 the inverse might have a pole on the line
 $\{\Re z=\frac{n+1}{2}-\gamma-\mu\}$;  we fix $\varepsilon_0>0$ such that
 there exists no pole with real part between $\frac{n+1}{2}-\gamma-\mu$ and 
 $\frac{n+1}{2}-\gamma-\mu\pm2\varepsilon_0$. In Propositions \ref{prop2} and
\ref{bl}, we next construct right and left parametrices modulo
parameter-dependent zero
order Green operators. Using conditions (E2) and (E3) we then can express 
$(\eta^\mu-\underline{ \widehat{A}})^{-1}$ as a principal edge symbol, 
and finally understand
$(\eta^\mu-\underline  A)^{-1}$. 
%% Theorem \ref{roughparametrix2} yields: 
 
 \begin{proposition}\label{prop2}
  For any $0<\varepsilon<\varepsilon_0$ there exists a
  parameter-dependent cone operator $b_R(\eta)\in
  C^{-\mu}(\Sigma;\gamma-\varepsilon,\gamma+\mu-\varepsilon,\mu)$ such that
   \begin{equation}\label{right}
     a(\eta)b_R(\eta)-1=g_R(\eta)\in 
     C^0_G(\Sigma;\gamma-\varepsilon,1-2\varepsilon,O;
     \gamma+\varepsilon,1-2\varepsilon,Q) 
   \end{equation}
 with a suitable asymptotic type 
 $Q\in\mbox{\rm As}(\gamma+\varepsilon,1-2\varepsilon)$. Moreover, for each
 $\eta$, the action of $b_R(\eta)$ on $\calH^{0,\gamma}_p(\bz)$ is independent
 of $\varepsilon$ and yields a bounded mapping 
  $$b_R(\eta):\calH^{0,\gamma}_p(\bz)\longrightarrow\calD(A_{\min}).$$
 \end{proposition}
 
 \begin{proof}
  %Let us first view $a(\eta)$ as an element in $C^\mu_O(\Sigma)$.
  According to Theorem \ref{roughparametrix1} there exists a flat
  parameter-dependent cone operator $\widetilde{b}(\eta)\in C^{-\mu}_O(\Sigma)$
  such that 
   \begin{equation}\label{formel1}
    a(\eta)\widetilde{b}(\eta)-1\equiv
    \omega(t[\eta])\,\op_M(\widetilde{f})(\eta)\,\omega_0(t[\eta])
    \qquad\text{mod }C^0_G(\Sigma)_\infty
   \end{equation}
  with a holomorphic Mellin symbol 
  $\widetilde{f}\in M^{-\infty}(\partial\bz;\Sigma)$. Setting
  $\widetilde{f}_0(z)=\widetilde{f}(z,0)$, we have  
   \begin{equation}\label{formel2}
    \omega(t[\eta])\left\{
    \op_M^{\widetilde{\gamma}-\frac{n}{2}}(\widetilde{f})(\eta)-
    \op_M^{\widetilde{\gamma}-\frac{n}{2}}(\widetilde{f}_0)\right\}
    \omega_0(t[\eta])\in
    C^0_G(\Sigma;\widetilde{\gamma},1,O;\widetilde{\gamma},1,O)
   \end{equation}
  for any $\widetilde{\gamma}\in\rz$. In fact, by a Taylor expansion, 
  $\widetilde{f}(z,t\eta)=\widetilde{f}_0(z)+\sum_{i=1}^n
  (t\eta_i)\widetilde{f}_i(z,t\eta)$ for suitable
  $\widetilde{f}_i(z,\eta)\in\ci(\Sigma,M^{-\infty}_O(\partial\bz))$. Therefore,
  the operator-family in \eqref{formel2} pointwise has the mapping properties
  of Green symbols from 
  $R^0_G(\Sigma;\widetilde{\gamma},1,O;\widetilde{\gamma},1,O)$,
  cf.\ Definition \ref{edgegreen}. Moreover, it is homogeneous of order 0 for
  large $|\eta|$ in the sense of \eqref{twisted} with respect to the group
  action of \eqref{groupaction}, hence also satisfies the required symbol
  estimates. 
  
  Next, we are going to modify $\widetilde{b}(\eta)$ by a smoothing Mellin term,
  i.e.\ we set 
   \begin{equation}\label{formel2.5}
     b_R(\eta)=\widetilde{b}(\eta)+
     \omega(t[\eta])\,t^\mu\,\op_M^{\gamma-\varepsilon-\frac{n}{2}}(f)\,
     \omega_0(t[\eta]),
   \end{equation}
  where we determine $f$ in such a way that the conormal symbol of
  $a(\eta)b_R(\eta)-1$ vanishes, i.e.\ 
   $$0=\sigma^0_M(ab_R)(z)-1=\sigma^0_M(a\widetilde{b})(z)+
     \sigma^\mu_M(a)(z-\mu)f(z)-1=\widetilde{f}_0(z)+
     \sigma^\mu_M(a)(z-\mu)f(z).$$
  Solving for $f$ yields 
   \begin{equation}\label{formel2.6}
    f(z)=-\sigma_M^\mu(a)^{-1}(z-\mu)\widetilde{f}_0(z)=
          \sigma_M^\mu(a)^{-1}(z-\mu)-\sigma_M^{-\mu}(\widetilde{b})(z). 
   \end{equation}
  By our choice of $\varepsilon_0$ in connection with Lemma \ref{shift3}, the 
  action of $b_R(\eta)$ on $\calH^{s,\gamma}_p(\bz)$ is independent
  of $\varepsilon$. Moreover, the description of $\calD(\amin)$
  given in Proposition \ref{mindomain} together with \eqref{formel2.6} implies
  that $b_R(\eta)$ maps $\calH^{0,\gamma}_p(\bz)$ to 
$\calD(A_{\min})$.
  Choose a cut-off  function $\omega_1$ with $\omega\omega_1=\omega_1$. 
  We claim that 
   \begin{equation}\label{formel3}
     a(\eta)b_R(\eta)-1\equiv
     \omega_1(t[\eta])\,\op_M^{\gamma-\varepsilon-\frac{n}{2}}
     (\widetilde{f}_0)\,\omega_0(t[\eta])+
     \omega_1(t[\eta])\,a(\eta)\,\omega(t[\eta])\,t^\mu\,
     \op_M^{\gamma-\varepsilon-\frac{n}{2}}(f)\,\omega_0(t[\eta])
   \end{equation}
  modulo $C^0_G(\Sigma;\gamma-\varepsilon,1-2\varepsilon,O;
  \gamma+\varepsilon,1-2\varepsilon,O)$. Indeed, this follows from
  \eqref{formel1} and \eqref{formel2} together with the fact that
  changing $\omega$ in \eqref{formel1} to $\omega_1$ only causes a remainder
  in $C^0_G(\Sigma)_\infty$, and  Lemma \ref{weight} (for 
  $\gamma-\varepsilon$ instead of $\gamma$). 

Since 
  $a(\eta)=\eta^\mu-t^{-\mu}\op_M(\sigma_M^\mu(a))$ on
  $[0,1]\times\partial\bz$, an application of \eqref{kommu} shows that the
second term on the right-hand side of
  \eqref{formel3} equals 
  \begin{align*}
   -&\omega_1(t[\eta])\,\op_M^{\gamma-\varepsilon-\frac{n}{2}}(\widetilde{f}_0)
     \,\omega_0(t[\eta])+\omega_1(t[\eta])\,(t\eta)^\mu\,
     \op_M^{\gamma-\varepsilon-\frac{n}{2}}(f)\,\omega_0(t[\eta])\\
   -&\omega_1(t[\eta])\,\op_M^{\gamma-\varepsilon-\frac{n}{2}}
     (T^{-\mu}\sigma_M^\mu(a))\,(1-\omega)(t[\eta])\,
     \op_M^{\gamma-\varepsilon-\frac{n}{2}}(f)\,\omega_0(t[\eta]).
  \end{align*}
  %Here, $T^{-\mu}$ is the operator of shifting by $-\mu$, cf.\ Theorem 
  %\ref{flatalgebra}. 
  The first term cancels with the first term of
  \eqref{formel3}, the other two belong to
  $C^0_G(\Sigma;\gamma-\varepsilon,1,O;\gamma+\varepsilon,1,Q)$, where $Q$ is
  the asymptotic type induced by the meromorphic structure of
  $(T^{-\mu}\sigma_M^\mu(a)^{-1})^*$. 
 \end{proof}
 
 In a similar way, one can also construct a rough left-parametrix for
 $a(\eta)$:

 \begin{proposition}\label{bl}
  Let $0<\varepsilon<\varepsilon_0$. Then there exists a
  parameter-dependent cone operator $b_L(\eta)\in
  C^{-\mu}(\Sigma;\gamma-\mu+\varepsilon,\gamma+\varepsilon,\mu)$ such that 
   \begin{equation}\label{left}
     b_L(\eta)a(\eta)-1=g_L(\eta)\in 
     C^0_G(\Sigma;\gamma+\varepsilon,\gamma+\varepsilon,\mu).
   \end{equation}
 \end{proposition}
 
 Passing in \eqref{right}, \eqref{left} to the principal edge symbols, cf.\
 \eqref{edgesymbol}, and solving for $(\eta^\mu-\underline{\widehat{A}})^{-1}$
 yields 
   \begin{equation}\label{formel4}
    (\eta^\mu-\underline{\widehat{A}})^{-1}=\sigma_\wedge^{-\mu}(b_R)(\eta)-
    \sigma_\wedge^{-\mu}(b_L)(\eta)\sigma_\wedge^0(g_R)(\eta)+
    \sigma_\wedge^0(g_L)(\eta)(\eta^\mu-\underline{\widehat{A}})^{-1}
    \sigma_\wedge^0(g_R)(\eta)
   \end{equation}
 on $\calK^{0,\gamma}_p(\partial\bz^\wedge)$. We are now going to show that
 the second and third term on the right-hand side of \eqref{formel4} are the
 principal edge symbol of a parameter-dependent Green operator. Let us set 
  \begin{equation}\label{eprime}
   \underline{\calE}^\prime=\underline{\calE}\oplus 
   \mbox{\rm im}\,\omega
   \left(\op_M^{\gamma+\mu-\varepsilon-\frac{n}{2}}(\sigma^\mu_M(A)^{-1})-
   \op_M^{\gamma+\mu+\varepsilon-\frac{n}{2}}(\sigma^\mu_M(A)^{-1})\right), 
  \end{equation}
 i.e.\ we add to $\underline{\calE}$ the asymptotic terms coming from the poles
 of $\sigma^\mu_M(A)^{-1}$ on the line $\re z=\frac{n+1}{2}-\gamma-\mu$. Note
 that these terms belong to $\calD(A_{\min})$. 
 
 \begin{proposition}\label{prop3}
  There exists a Green symbol 
  $g\in R^{-\mu}_G(\Sigma;\gamma,\varepsilon;\gamma,\mu+\varepsilon)$ for some 
  $\varepsilon>0$ such that 
   $$\sigma_\wedge^{-\mu}(g)(\eta)=
     -\sigma_\wedge^{-\mu}(b_L)(\eta)\sigma_\wedge^0(g_R)(\eta)+
     \sigma_\wedge^0(g_L)(\eta)(\eta^\mu-\underline{\widehat{A}})^{-1}
     \sigma_\wedge^0(g_R)(\eta)$$
  and the integral kernel $k_g$ of $g$, cf. \mbox{\rm\eqref{kernel}}
  and Theorem \mbox{\rm\ref{characterization}}, satisfies 
   $$k_g(\eta,t,x,t^\prime,x^\prime)\in 
     S^{-\mu}_{cl}(\Sigma,
      \calS^\gamma_{\underline{\calE}^\prime}(\partial\bz^\wedge)\pit
      \calS^{-\gamma}_{\varepsilon/2}(\partial\bz^\wedge))=
     S^{-\mu}_{cl}(\Sigma)\pit
      \calS^\gamma_{\underline{\calE}^\prime}(\partial\bz^\wedge)\pit
      \calS^{-\gamma}_{\varepsilon/2}(\partial\bz^\wedge),$$
  where we set $\calS^\gamma_{\underline{\calE}^\prime}(\partial\bz^\wedge)=
  \calS^{\gamma+\mu}_{\varepsilon/2}(\partial\bz^\wedge)\oplus
  \underline{\calE}^\prime$. 
%and $\underline{\calE}^\prime$ is as described in
%  \eqref{eprime}.
 \end{proposition}
 
 \begin{proof}
  For brevity let us write $\widehat{g}(\eta)=\sigma_\wedge^{-\mu}(g)(\eta)$,
  and analogously for the principal edge symbols of $b$, $g_L$, and $g_R$. 
  By the previous proposition, we have 
   $$\calK^{s,\gamma}_2(\partial\bz^\wedge)\hookrightarrow
     \calK^{s,\gamma-\varepsilon}_2(\partial\bz^\wedge)
     \xrightarrow{\widehat{g}_R(\eta)}
     \calS^{\gamma+\varepsilon}_Q(\partial\bz^\wedge)
     \hookrightarrow
     \calS^{\gamma-\mu+\varepsilon}_Q(\partial\bz^\wedge),$$
  where one considers 
  $Q\in\mbox{\rm As}(\gamma+\varepsilon,1-2\varepsilon)$ as an asymptotic type
  $Q\in\mbox{\rm As}(\gamma-\mu+\varepsilon,\mu)$ for the last embedding. 
By standard mapping 
  properties of cone operators there exists an asymptotic
  type $Q^\prime\in\mbox{\rm As}(\gamma+\varepsilon,\mu)$ such that
  $\widehat{b}_L(\eta):\calS^{\gamma-\mu+\varepsilon}_Q(\partial\bz^\wedge)\to 
  \calS^{\gamma+\varepsilon}_{Q^\prime}(\partial\bz^\wedge)\hookrightarrow
  \calS^{\gamma}_{Q^\prime}(\partial\bz^\wedge)$, considering 
  $Q^\prime$ as an element of  
  $\mbox{\rm As}(\gamma,\mu+\varepsilon)$. 
Making similar considerations for
  the adjoint, we thus obtain  
  $\widehat{b}_L(\eta)\widehat{g}_R(\eta)\in 
  R^{(-\mu)}_G(\Sigma;\gamma,\varepsilon;\gamma,\mu+\varepsilon)$.
    
  Next, $(\eta^\mu-\underline{\widehat{A}})^{-1}$ is a smooth function on
  $\Sigma\setminus\{0\}$ with values in 
  $\calL(\calK^{0,\gamma}_p(\partial\bz^\wedge),
  \calD(\underline{\widehat{A}}))$. In view of assumption \mbox{\rm(E2)} on the
  scaling invariance of $\underline{\calE}$, it is homogeneous of order $-\mu$
  in the sense of \eqref{twisted}. In particular, 
   $(\eta^\mu-\underline{\widehat{A}})^{-1}\in 
     S^{(-\mu)}(\Sigma;\calK^{0,\gamma}_p(\partial\bz^\wedge),
     \calK^{0,\gamma+\varepsilon}_p(\partial\bz^\wedge)$
  for sufficiently small $\varepsilon>0$. But then it is clear that 
  $\widehat{g}_L(\eta)(\eta^\mu-\underline{\widehat{A}})^{-1}
  \widehat{g}_R(\eta)$ also belongs to 
  $R^{(-\mu)}_G(\Sigma;\gamma,\varepsilon;\gamma,\mu+\varepsilon)$. If we now 
  define 
   $$g(\eta)=\chi(\eta)\{-\widehat{b}_L(\eta)\widehat{g}_R(\eta)+
     \widehat{g}_L(\eta)(\eta^\mu-\underline{\widehat{A}})^{-1}
     \widehat{g}_R(\eta)\}$$ 
  with an arbitrary zero excision function $\chi(\eta)$, then $g\in  
  R^{-\mu}_G(\Sigma;\gamma,\varepsilon;\gamma,\mu+\varepsilon)$, and
  the principal edge symbol is given by the formula stated in the proposition.
  It remains to investigate the kernel. 
  
  Since $g(\eta)$ is a Green symbol of the given class (and by the kernel
  characterization) there exists some asymptotic type $Q\in\mbox{\rm
  As}(\gamma,\mu+\frac{\varepsilon}{2})$ such that the integral kernel $k_g$ of
  $g(\eta)$ belongs to $S^{-\mu}_{cl}(\Sigma)\pit
  \calS^{\gamma+\varepsilon}_Q(\partial\bz^\wedge)\pit
  \calS^{-\gamma}_{\varepsilon/2}(\partial\bz^\wedge)$. According to Definition
  \ref{subspaces}, 
   $$\calS^{\gamma+\varepsilon}_Q(\partial\bz^\wedge)=
     \calS^{\gamma+\varepsilon/2}_\mu(\partial\bz^\wedge)\oplus
     \calE_Q.$$
  By possibly shrinking $\varepsilon$, we may assume that $Q$ contains no
  triple $(q,l,L)$ with $\re q<\frac{n+1}{2}-\gamma-\mu$. By possibly
  enlarging $Q$ we can assume that $\underline{\calE}^\prime\subset
  \calE_Q$, and therefore $\calE_Q=\underline{\calE}^\prime\oplus V$ for a 
  certain finite-dimensional space $V$. Therefore $k_g$ can be written as
  $k_g=k_g^0+k_g^1$ with $k_g^0$ containing
  the contribution of $\underline{\calE}^\prime$, and $k_g^1$ that of $V$. 
  However, from identity \eqref{formel4} one sees that $g(\eta)$ maps into the
  domain of $\underline{\widehat{A}}$, and 
  therefore $k_g^1$ must equal $0$. This then finishes the proof, since 
  $\calS^{\gamma+\varepsilon/2}_\mu(\partial\bz^\wedge)=
  \calS^{\gamma+\mu}_{\varepsilon/2}(\partial\bz^\wedge)$. 
 \end{proof}
 
 With $b_R(\eta)$ from Proposition \ref{prop2} and $g(\eta)$ from Proposition 
 \ref{prop3} let us now define 
  \begin{equation}\label{formel4.5}
    b(\eta)=b_R(\eta)+\sigma\,g(\eta)\,\sigma_0=\widetilde{b}(\eta)+
    \omega(t[\eta])\,t^\mu\,\op_M^{\gamma-\varepsilon-\frac{n}{2}}(f)\,
    \omega_0(t[\eta])+\sigma\,g(\eta)\,\sigma_0
  \end{equation}
 wit cut-off functions $\sigma,\sigma_0$. Clearly, $b(\eta): 
\calH^{0,\gamma}_p(\bz)\to\calD(\underline{A})$ 
 for each $\eta$. As explained
 in Definition \ref{flatcone} and the subsequent paragraph, there are then
 $h$, $p$, and a cut-off function $\sigma_1$ such that 
  $$b(\eta)=\sigma
    \left\{t^\mu\,\op_M^{\gamma-\frac{n}{2}}(h)(\eta)+
    \omega(t[\eta])\,t^\mu\,\op_M^{\gamma-\varepsilon-\frac{n}{2}}(f)\,
    \omega_0(t[\eta])+g(\eta)\right\}\sigma_0+
    (1-\sigma)\,p(\eta)\,(1-\sigma_1).$$
 Since $A$ has $t$-independent coefficients on $[0,1]\times\partial\bz$,
 we can also choose $h$ to be independent of $t$, i.e.\ $h\in
 M^{-\mu}_O(\partial\bz;\Sigma)$. Then 
\begin{equation} \label{br} 
\sigma_\wedge^{-\mu}(b_R)(\eta)=
    t^\mu\,\op_M^{\gamma-\frac{n}{2}}(h)(\eta)+
    \omega(t|\eta|)\,t^\mu\,\op_M^{\gamma-\varepsilon-\frac{n}{2}}(f)\,
    \omega_0(t|\eta|).
\end{equation}
 As $\sigma_\wedge^{-\mu}(g)(\eta)=g(\eta)$ for large $|\eta|$ by
 construction of, we deduce from \eqref{formel4} and  \eqref{br}  that  
  \begin{equation}\label{formel5}
   b(\eta)=\sigma\,(\eta^\mu-\underline{\widehat{A}})^{-1}\,\sigma_0+
   (1-\sigma)\,p(\eta)\,(1-\sigma_1),\qquad|\eta|\ge R, 
  \end{equation}
 for a sufficiently large $R>0$. Changing the cut-off functions 
 only alters $b(\eta)$ by a nice remainder: 
 
 \begin{lemma}\label{lemma1}
  Let $\widetilde{\sigma},\widetilde{\sigma}_0,\widetilde{\sigma}_1$ be
  cut-off functions satisfying the conditions posed in Definition
  \mbox{\rm\ref{flatcone}}  and the subsequent paragraph. Then 
   $$b(\eta)=\widetilde{\sigma}\,(\eta^\mu-\underline{\widehat{A}})^{-1}\,
     \widetilde{\sigma}_0+
     (1-\widetilde{\sigma})\,p(\eta)\,(1-\widetilde{\sigma}_1)+
     r(\eta),\qquad|\eta|\ge R,$$ 
  with a remainder $r(\eta)$ having an integral kernel 
  $($with a certain $\varepsilon>0)$
   $$k_r(\eta,y,y^\prime)\in\calS(\Sigma,
     \calC^{\infty,\gamma}_{\underline{\calE}^\prime}(\bz)\pit
     \calC^{\infty,-\gamma}_{\varepsilon}(\bz)).$$
  Here, we have set $\calC^{\infty,\gamma}_{\underline{\calE}^\prime}(\bz)=
  \calC^{\infty,\gamma+\mu}_{\varepsilon}\oplus\underline{\calE}^\prime$. 
  Hence we have not only $r(\eta)\in
  C^{-\infty}_G(\Sigma;\gamma,\gamma,\varepsilon)$, but also 
  $r(\eta)\in\calS(\Sigma,\calL(\calH^{0,\gamma}_p(\bz),
  \calD(\underline{A})))$.
 \end{lemma}
 
 \begin{proof}
  This statement is easily seen using the representation of $b(\eta)$ in
  \eqref{formel4.5}. Changing cut-off functions alters $\widetilde{b}(\eta)$ 
  only by a flat Green symbol in $C^{-\infty}_G(\Sigma)_\infty$, which has in
  particular a kernel of the mentioned structure. It remains to note that
  $\sigma g(\eta)\sigma_0-\widetilde{\sigma}g(\eta)\widetilde{\sigma}_0=
  \sigma g(\eta)(\sigma_0-\widetilde{\sigma}_0)+
  (\widetilde{\sigma}-\sigma)g(\eta)\widetilde{\sigma}_0$ and both these terms 
  have the required structure (recall that they are Green symbols of order
  $-\infty$, since both $\widetilde{\sigma}-\sigma$ and
  $\sigma_0-\widetilde{\sigma}_0$ belong
  to $\cicomp(]0,1[)$). 
 \end{proof}
 
 \begin{proposition}\label{prop4} 
  If $b(\eta)$ is as in \eqref{formel5}, then for $|\eta|$ large enough, 
   $$(\eta^\mu-\underline{A})b(\eta)-1=r_R(\eta),\qquad 
     b(\eta)(\eta^\mu-\underline{A})-1=r_L(\eta),$$
  with 
   $$r_R(\eta)\in C^{-\infty}_G(\Sigma;\gamma,\gamma,\varepsilon),\qquad 
     r_L(\eta)\in\calS(\Sigma,\calL(\calD(\underline{A})))$$ 
  for suitable $\varepsilon>0$. In particular,
  $(\eta^\mu-\underline{A}):\calD(\underline{A})\to\calH^{0,\gamma}_p(\bz)$ is
  invertible for  sufficiently large $|\eta|$. 
 \end{proposition}
 
 \begin{proof}
  Since $A$ is a local operator, we can write 
  $(\eta^\mu-\underline{A})=
  \sigma(\eta^\mu-\underline{A})\sigma_0+
  (1-\sigma)(\eta^\mu-A)(1-\sigma_1)$ 
  with cut-off functions as above. Then 
   $$(\eta^\mu-\underline{A})\,b(\eta)=
   \sigma\,(\eta^\mu-\underline{A})\,\sigma_0\,b(\eta)+
     (1-\sigma)\,(\eta^\mu-\underline{A})\,(1-\sigma_1)\,b(\eta).$$ 
  To treat the first summand choose a representation of $b(\eta)$ as in Lemma
  \ref{lemma1} with $\widetilde{\sigma}\sigma_0=\sigma_0$. Then 
   $$\sigma\,(\eta^\mu-\underline{A})\,\sigma_0\,b(\eta)=\sigma+
     \sigma\,(\eta^\mu-\underline{A})\,\sigma_0\,r(\eta),$$
  where the second term on the right-hand side belongs to
  $C^{-\infty}_G(\Sigma;\gamma,\gamma,\varepsilon)$, since $A$ maps
  $\calC^{\infty,\gamma}_{\underline{\calE}^\prime}$ to $\cicomp(\intb)$. 
  For the second summand we choose for $b(\eta)$ a representation with
  $\widetilde{\sigma}_0$ such that 
  $\sigma_1\widetilde{\sigma}_0=\widetilde{\sigma}_0$. Then 
   $$(1-\sigma)\,(\eta^\mu-\underline{A})\,(1-\sigma_1)\,b(\eta)\equiv
     (1-\sigma)+(1-\sigma)\,((\eta^\mu-A)p(\eta)-1)\,(1-\widetilde{\sigma}_1)$$
  modulo a remainder of the prescribed form. Since $p(\eta)$ is a
parametrix to $\eta^\mu-\underline
A$ in the interior of $\bz$,
the second term on the right-hand side is of the form
$(1-\sigma)d(\eta)(1-\widetilde{\sigma}_1)$
  with $d(\eta)\in L^{-\infty}(2\bz;\Sigma)$, hence is a
  remainder of the desired form. The considerations for
  $b(\eta)(\eta^\mu-\underline{A})$ are analogous. 
 \end{proof}

 To finish the proof of Theorem \ref{maintheorem}.b), it remains to modify 
 $b(\eta)$ in such a way that we obtain the inverse of 
 $\eta^\mu-\underline{A}$. To do so, we may assume that $r_R(\eta)$ of 
 Proposition \ref{prop4} satisfies 
 $\|r_R(\eta)\|_{\calL(\calH^{0,\gamma}_p(\bz))}\le\frac{1}{2}$ for all 
 $\eta\in\Sigma$ (otherwise we multiply $r_R(\eta)$ with a suitable zero 
 excision function $\chi(\eta)$). Then 
 $1+r_R(\eta)$ is invertible in $\calH^{0,\gamma}_p(\bz))$ for all 
 $\eta\in\Sigma$ and 
  $$(1+r_R(\eta))^{-1}=1-r_R(\eta)+r_R(\eta)(1+r_R(\eta))^{-1}r_R(\eta)=:
    1+s(\eta).$$
 Clearly, $s(\eta)$ belongs to  
 $C^{-\infty}_G(\Sigma;\gamma,\gamma,\varepsilon)$, again. Hence, by 
 Proposition \ref{prop4} and \eqref{formel4.5},
  \begin{equation}\label{efgh}
   (\eta^\mu-\underline{A})^{-1}=b(\eta)(1+r_R(\eta))^{-1}=
   b(\eta)+b(\eta)s(\eta)=\widetilde{b}(\eta)+\widetilde{r}(\eta)
  \end{equation}
 for large $|\eta|$, with $\widetilde{r}(\eta)\in
 C^{-\mu}_G(\Sigma;\gamma,\gamma,\varepsilon)$ 
 (note that the smoothing Mellin term in \eqref{formel4.5} belongs to 
 $C^{-\mu}_{M+G}(\Sigma;\gamma-\varepsilon,\gamma+\mu-\varepsilon,\mu)
 \subset C^{-\mu}_G(\Sigma;\gamma,\gamma,\varepsilon)$). 

%%%%%%%%%%%%%%%%%%%%%%%%%%%%%%%%%%%%%%%%%%%%%%%%%%%%%%%%%%%%%%%%%%%%%
%%%%%%%%%%%%%%%%%%%%%%%%%%%%%%%%%%%%%%%%%%%%%%%%%%%%%%%%%%%%%%%%%%%%%

\section{Spectral invariance and bounded imaginary powers}\label{section3}

\subsection{Independence on the choice of $1<p<\infty$}\label{section3.1}

A closed extension $\underline A$ of $A$ in
$\calH^{0,\gamma}_{p_0}(\bz)$
with domain $\calD(\amin)\oplus \underline \calE$ as in  
Theorem \ref{maxdomain} induces
closed extensions $\underline A_p$ in
$\calH^{0,\gamma}_p(\bz)$ for all $1<p<\infty$ by letting $A$
act on  
% Let $\underline{A}_p$ deote a closed extension of $A$ in
% $\calH^{0,\gamma}_p(\bz)$ defined by the domain
  $$\calD(\underline{A}_p)=\calD_p(A_{\min})\oplus\underline{\calE},$$
% cf.\ Theorem \ref{maxdomain}, 
where
 $\calD_p(A_{\min})$ is the domain of the closure of $A$ in
$\calH^{0,\gamma}_p(\bz)$.
 
 \begin{theorem}\label{independence1}
  Suppose $\underline{A}=\underline{A}_{p_0}$ satisfies the assumptions of
  Theorem {\rm\ref{maintheorem}} for one $p_0$, and let $c(\eta)$ be the
  inverse of $\eta^\mu-\underline{A}_{p_0}$ for $|\eta|\ge R$. Then $c(\eta)$,
  $|\eta|\ge R$, also yields the inverse of $\eta^\mu-\underline{A}_p$ for
  all $1<p<\infty$. 
 \end{theorem}
 \begin{proof}
  By \eqref{efgh}, $c(\eta)=b(\eta)(1+r_R(\eta))^{-1}$ with
  $b(\eta)=b_R(\eta)+\sigma_0\,g(\eta)\,\sigma_0$. According to Proposition 
  \ref{prop2}, $b_R(\eta)$ maps $\calH^{0,\gamma}_p(\bz)$ into
  $\calD_p(A_{\min})$. By Proposition \ref{prop3},
  $\sigma_0\,g(\eta)\,\sigma_0$
  maps into $\calC^{\infty,\gamma}_{\underline{\calE}^\prime}(\bz)$,
  which is a subspace of
  $\calD(\underline{A}_p)$. Hence 
   $$c(\eta):\calH^{0,\gamma}_p(\bz)\longrightarrow
     \calD(\underline{A})=\calD(\underline{A}_p)$$
  for all $1<p<\infty$. Moreover, $c(\eta)(\eta^\mu-\underline{A}_p)=1$ on 
  $\cicomp(\intb)\oplus\underline{\calE}$ for each $p$, since this is true for
  $p=p_0$ and the left-hand side of the latter equation is independent of $p$
  on $\cicomp(\intb)\oplus\underline{\calE}$. Similarly,  
  $(\eta^\mu-\underline{A}_p)c(\eta)=1$ on $\cicomp(\intb)$ for all $p$. Thus a
  density argument gives the result. 
 \end{proof}

 Equation \eqref{formel4} together with Proposition \ref{prop3}
 shows that 
  $(\eta^\mu-\underline{\widehat{A}}_p)^{-1}$ is the sum of two principal edge
symbols, namely those of an operator in $C^{-\mu}_O(\Sigma)$ and one in 
$R^{-\mu}_G(\Sigma;\gamma,\gamma,\varepsilon)$. Arguing similarly as above, we
therefore obtain

%\in 
%    R^{(-\mu)}(\Sigma)+
%    R^{(-\mu)}_G(\Sigma;\gamma,\gamma,\varepsilon),$$
% where $R^{(-\mu)}(\Sigma)$ denotes the space of all principal edge symbols
% $\sigma^{-\mu}_\wedge(c)(\eta)$ with $c\in C^{-\mu}_O(\Sigma)$, and 
% $R^{(-\mu)}_G(\Sigma;\gamma,\gamma,\varepsilon)$ the space of all homogeneous
% principal symbols $g_{(-\mu)}(\eta)$ with  
% $g\in R^{-\mu}_G(\Sigma;\gamma,\gamma,\varepsilon)$. Then, arguing
% similarly as above, we obtain
 
 \begin{theorem}\label{independence2}
  If $\underline{A}=\underline{A}_p$ satisfies ellipticity assumptions 
  \mbox{\rm(E1)}, \mbox{\rm(E2)}, and \mbox{\rm(E3)} of Section
  \mbox{\rm\ref{section2.1}} for one $p$, then automatically for all
  $1<p<\infty$. 
 \end{theorem}

%%%%%%%%%%%%%%%%%%%%%%%%%%%%%%%%%%%%%%%%%%%%%%%%%%%%%%%%%%%%%%%%%%%%%%

\subsection{Bounded imaginary powers}\label{section3.2}
 
 In the paper \cite{CSS1} we have shown that the closure of a cone 
 differential operator -- under ellipticity conditions \mbox{\rm(E1)} and 
 \mbox{\rm(E3)} with $\underline{\calE}=\{0\}$ -- posseses bounded imaginary 
 powers whose operator-norm in $\calH^{0,\gamma}_p(\bz)$ can be estimated by 
 $c_p e^{\theta|z|}$, where $\theta$ is the angle determining 
 $\Lambda=\Lambda_\theta$. We also had pointed out in Remark 5.5 of \cite{CSS1}
 that the validity of this result `only' relies on the structure of the 
 resolvent and not on the fact that we dealt with the minimal extension.  
  Theorem \ref{maintheorem} now states that the resolvent of a 
 general closed extension has exactly this required structure (in 
 \cite{CSS1} we described the resolvent $(\lambda-\overline{A})^{-1}$ in 
 terms of anisotropic symbols, while here we described 
 $(\eta^\mu-\underline{A})^{-1}$; both ways, however, are obviously 
 equivalent). Thus we have the following result: 
 
 \begin{theorem}\label{powers}
  Let $\underline{A}$ be a closed extension of a cone differential  operator 
  $A$, satisfying the ellipticity assumptions \mbox{\rm(E1)}, \mbox{\rm(E2)}, 
  and \mbox{\rm(E3)} with respect to $\Lambda=\Lambda_\theta$.
  Then the resolvent $(\lambda- \underline A)^{-1}$ exists for large $\lambda$ 
  in $\Lambda$, and its norm in $\calH^{0,\gamma}_p(\bz)$ decays like
  $1/|\lambda|$. Moreover, there exists a constant  $c\ge 0$ such that 
  $\underline{A}+c$ has bounded imaginary powers and, 
  for a certain constant $c_p\ge0$, 
   $$\|(\underline{A}+c)^{i\varrho}\|_{\calL(\calH^{0,\gamma}_p(\bz))}\le 
     c_p\,e^{\theta|\varrho|}\qquad\forall\;\varrho\in\rz.$$ 
  As the construction of complex powers shows, we can take $c=0$ if
  $\underline{A}$ has no spectrum in $\Lambda_\theta$. 
 \end{theorem} 

 Let us mention that the operator $A+c$ does not satisfy the assumption of
 constant coefficients near the boundary (since we have to write
 $c=t^{-\mu}(t^\mu c)$). However, the structure of the resolvent remains
 uneffected by the shift with a constant $c$. 

%%%%%%%%%%%%%%%%%%%%%%%%%%%%%%%%%%%%%%%%%%%%%%%%%%%%%%%%%%%%%%%%%%%%%
%%%%%%%%%%%%%%%%%%%%%%%%%%%%%%%%%%%%%%%%%%%%%%%%%%%%%%%%%%%%%%%%%%%%%

\section{The Laplace-Beltrami operator}\label{section4}
 
 Let the interior of $\bz$ be equipped with a metric that coincides 
 with $dt^2+t^2g$  on  ${]0,1[}\times\partial\bz$ for some fixed 
 metric $g$ on $\partial\bz$ (straight conical degeneracy).
The associated Laplacian $\Delta$ is a second order cone differential 
operator, and 
  $$\Delta=t^{-2}\,\{(t\partial_t)^2+(n-1)t\partial_t+\Delta_\partial\}, 
    \qquad n=\text{\rm dim}\,\partial\bz,$$
 near the boundary of $\bz$, where $\Delta_\partial$ denotes the Laplacian 
 on $\partial\bz$ with respect to $g$. 
 
 Clearly, $-\Delta$ satisfies ellipticity condition {\rm(E1)} of Section 
 \ref{section2.1} for any sector $\Lambda$ not containing positive reals. 
 
%%%%%%%%%%%%%%%%%%%%%%%%%%%%%%%%%%%%%%%%%%%%%%%%%%%%%%%%%%%%%%%%%%%%%

\subsection{The conormal symbol}\label{section4.1}

 %By the  definition given in the beginning of Section \ref{section1.1}, the 
 %conormal symbol of the Laplacian is 
Let us first analyze  the inverse of the conormal symbol 
  $$\sigma^2_M(\Delta)(z)=z^2-(n-1)z+\Delta_\partial.$$
%Let us investigate the inverse of this function. 
To this end denote by
 $0=\lambda_0>\lambda_1>\ldots$ the eigenvalues of $\Delta_\partial$ and by
 $E_0,\,E_1,\ldots$ the corresponding eigenspaces. Moreover, let 
 $\pi_j\in\calL(L_2(\partial\bz))$ be the orthogonal projection onto 
 $E_j$. 
 
 The {\em non}-bijectivity points of $\sigma^2_M(\Delta)$ are exactly the 
 points $z=q_j^+$ and $z=q_j^-$ with 
  \begin{equation}\label{pjpm}
   q_j^\pm=\mbox{$\frac{n-1}{2}\pm
   \sqrt{\big(\frac{n-1}{2}\big)^2-\lambda_j}$},
   \qquad j\in\nz_0.
  \end{equation}
 Note the symmetry $q_j^+=(n-1)-q_j^-$. It is straightforward to 
 calculate that 
  $$(z^2-(n-1)z+\Delta_\partial)^{-1}=
    \smsum_{j=0}^\infty\mbox{$\frac{1}{q_j^+-q_j^-}$}\Big(
    \mbox{$\frac{1}{z-q_j^+}$}-\mbox{$\frac{1}{z-q_j^-}$}\Big)\pi_j.$$
 Hence, in case $\text{\rm dim}\,\bz\ge3$,
  $$(z^2-(n-1)z+\Delta_\partial)^{-1}\equiv
    \mbox{$\pm\frac{1}{q_j^+-q_j^-}$}\,\pi_j\,(z-q_j^\pm)^{-1}\qquad
    \text{near }z=q_j^\pm$$
 modulo holomorphic germs. In case $\text{\rm dim}\,\bz=2$ the same formula
 holds near $z=q_j^\pm$ if $j\ge1$ but, since then $q_0^+=q_0^-=0$,  
  $$(z^2+\Delta_\partial)^{-1}\equiv\pi_0\,z^{-2}
  \qquad\text{near }z=0.$$
For $\text{\rm dim}\,\bz=1$, this simplifies to 
$\sigma^2_M(\Delta)^{-1}(z)=(z^2+z)^{-1}=\frac1z-\frac1{z+1}$ with only two
poles in $q^-_0 = -1$ and $q^+_0=0$ and associated `eigenspace' $E_0=\cz$.

%%%%%%%%%%%%%%%%%%%%%%%%%%%%%%%%%%%%%%%%%%%%%%%%%%%%%%%%%%%%%%%%%%%%%

\subsection{Maximal domain and dilation invariance}\label{section4.2}
 
 With $q_j^\pm$ we associate the function space 
  $$\calE_{q_j^\pm}=E_j\otimes\omega\,t^{-q_j^\pm}=
    \{e(x)\,\omega(t)\,t^{-q_j^\pm}\st e\in E_j\},\qquad j\in\nz,$$
 and for $q_0^\pm$ we set 
  \begin{equation}\label{ep0pm}
   \calE_{q_0^\pm}=
    \begin{cases}
     E_0\otimes\omega+E_0\otimes\omega\log t & 
      \text{\rm dim}\,\bz=2\\
     E_0\otimes\omega\,t^{q_0^\pm}& \text{\rm dim}\,\bz\not=2
    \end{cases}.
  \end{equation}
 Furthermore, for $\gamma\in\rz$, set 
  $$I_\gamma=\{q_j^\pm\st j\in\nz_0\}\cap
    \,\mbox{$]\frac{n+1}{2}-\gamma-2,\frac{n+1}{2}-\gamma[$}=
    \{q_j^\pm\st j\in\nz_0\}\cap
    \,\mbox{$]\frac{n-1}{2}-\gamma-1,\frac{n-1}{2}-\gamma+1[$}.$$
 Applying Theorems \ref{mindomain} and \ref{maxdomain} to
 $A=\Delta$, we get the following: 
 
 \begin{proposition}
  The domain of the maximal extension of $\Delta$ in
  $\calH^{0,\gamma}_p(\bz)$ is 
   $$\calD(\Delta_{\max})=\calD(\Delta_{\min})\oplus
     \mathop{\mbox{$\oplus$}}_{q\in I_\gamma}\calE_q.$$
  In case $q_j^\pm\not=\frac{n+1}{2}-\gamma-2$ for all $j$, the minimal 
  domain is $\calD(\Delta_{\min})=\calH^{2,\gamma+2}_p(\bz)$. 
 \end{proposition}

 Let us now describe the closed extensions $\underline{\Delta}$ of $\Delta$ 
 satisfying condition {\rm(E2)} of Section \ref{section2.1}. For convenience 
 we shall call such  extensions dilation invariant. A straightforward
 calculation (or an application of Lemmas 5.11 and 5.12 of \cite{GiMe}) 
 yields: 
 
 \begin{proposition}\label{invdomain}
  Consider $\Delta$ as an unbounded operator in $\calH^{0,\gamma}_p(\bz)$. 
  The dilation invariant extensions $\underline{\Delta}$ are precisely those 
  with a domain of the form 
   \begin{equation}\label{dilinv} 
    \calD(\underline{\Delta})=\calD(\Delta_{\min})\oplus
    \mathop{\mbox{$\oplus$}}_{q\in I_\gamma}\underline{\calE}_q,\qquad
    \underline{\calE}_q\text{ subspace of }\calE_q, 
   \end{equation}
  where in case $\mbox{\rm dim}\,\bz=2$ either 
  $\underline{\calE}_0=\{0\}$ or $\underline{\calE}_0=E_0\otimes\omega$ or 
  $\underline{\calE}_0=\calE_0$, cf.\ \eqref{ep0pm}. 
 \end{proposition}

 Let us point out that in \eqref{dilinv} the sum is taken over all  
 $q\in I_\gamma$ and that the summand $\underline{\calE}_q=\{0\}$ may occur
 several times. 
  
%%%%%%%%%%%%%%%%%%%%%%%%%%%%%%%%%%%%%%%%%%%%%%%%%%%%%%%%%%%%%%%%%%%%%

\subsection{Adjoint operators}\label{section4.3}
 
 Since the scalar-product $\skp{\cdot}{\cdot}_{0,0}$ of 
 $\calH^{0,0}_2(\bz)$ yields an identification of the dual space of 
 $\calH^{0,\gamma}_p(\bz)$ with $\calH^{0,-\gamma}_{p^\prime}(\bz)$, 
 the adjoint $\underline{\Delta}^*$ of an extension $\underline{\Delta}$ in 
 $\calH^{0,\gamma}_p(\bz)$ is the an unbounded operator in 
 $\calH^{0,-\gamma}_{p^\prime}(\bz)$ given by the action of 
 $\Delta$ on the domain 
  $$\calD(\underline{\Delta}^*)=
    \{v\in\calH^{0,-\gamma}_{p^\prime}(\bz)\st 
    \exists\,f\in\calH^{0,-\gamma}_{p^\prime}(\bz)\;
    \forall\,u\in\calD(\underline{\Delta}):\quad
    \skp{v}{\Delta u}_{0,0}=\skp{f}{u}_{0,0}\}.$$
 It is easy to see that $\Delta_{\min}^*=\Delta_{\max}$ and 
 $\Delta_{\max}^*=\Delta_{\min}$. 
 
 We shall now compute explicitly the adjoints of the dilation invariant
 extensions.  For an analysis of 
 adjoints of general cone differential operators (in case $p=2$) see  
 the paper \cite{GiMe}. Define  
  $$\skpa{\cdot}{\cdot}:
    \calD^\gamma_p(\Delta_{\max})\times
    \calD^{-\gamma}_{p^\prime}(\Delta_{\max})
    \longrightarrow\cz,\quad
    \skpa{u}{v}=\skp{\Delta u}{v}_{0,0}-\skp{u}{\Delta v}_{0,0},$$
 where the indices $\sigma,r$ in $\calD^\sigma_r$ now indicate that we 
 consider the Laplacian in the Sobolev space $\calH^{0,\sigma}_r(\bz)$. 
 Then the domain of the adjoint operator $\underline{\Delta}^*$ is just the 
 orthogonal space (with respect to this pairing) to the domain of 
 $\underline{\Delta}$, i.e. 
  $$\calD^{-\gamma}_{p^\prime}(\underline{\Delta}^*)=
    \calD^\gamma_p(\underline{\Delta})^\perp.$$ 
 Since $\skpa{u}{v}=0$ whenever $u$ or $v$ belong to the minimal domain, we 
 classify first which elements of $\oplus_{q\in I_{-\gamma}}\calE_q$ are
 orthogonal to a given element of $\oplus_{q\in I_{\gamma}}\calE_q$. 
 
 Let $u=e\,\omega\,t^{-q}$ with $q=q_j^+$ or $q=q_j^-$ and $e\in E_j$ for some
 fixed $j\in\nz_0$. If $v_\pm=f\,\omega\,t^{-q_k^\pm}$ with $f\in E_k$, an 
 elementary calculation yields 
  $$(\Delta u)\overline{v_\pm}-u\,\overline{\Delta v_\pm}= 
    2(q_k^\pm-q)\,e\,\overline{f}\,\omega\omega^\prime t^{-q-q_k^\pm-1},$$
 hence $\skpa{u}{v_\pm}=0$ if and only if $q_k^\pm=q$ or 
 $\skp{e}{f}_{L_2(\partial\bz)}=0$. 
 
 If $\mbox{\rm dim}\,\bz=2$ and $u=c\omega+d\omega\log t$ 
 with $c,d\in\cz$ and $v_\pm=f\,\omega\,t^{-q_k^\pm}$ with $f\in E_k$ and 
 $k\not=0$ then 
  $$(\Delta u)\overline{v_\pm}-u\,\overline{\Delta v_\pm}= 
    2\overline{f}(d+cq_k^\pm+dq_k^\pm\log t)\omega\omega^\prime 
    t^{-q_k^\pm-1},$$
 hence $\skpa{u}{v_\pm}=0$, since $\skp{1}{f}_{L_2(\partial\bz)}=0$. If $u$ 
 is as before and $v=c_0\omega+d_0\omega\log t$ with $c_0,d_0\in\cz$, then 
  $$(\Delta u)\overline{v}-u\,\overline{\Delta v}= 
    2(\overline{c}_0d-\overline{d}_0c)\omega\omega^\prime t^{-3},$$
 hence if both $u$ and $v$ are different from zero, $\skpa{u}{v}=0$ if  and 
 only if $\overline{v}$ is a multiple of $u$. 
 
 From this we derive the following description of adjoints of 
 dilation invariant extensions: 
 
 \begin{theorem}\label{adjoint}
  Let $\underline{\Delta}$ be a dilation invariant extension of $\Delta$ in 
  $\calH^{0,\gamma}_p(\bz)$ with domain 
   $$\calD^\gamma_p(\underline{\Delta})=
     \calD^\gamma_p(\Delta_{\min})\oplus
     \mathop{\mbox{$\oplus$}}_{q\in I_\gamma}\underline{\calE}_q$$
  as described in Proposition \mbox{\rm\ref{invdomain}}. Then the domain of 
  the adjoint $\underline{\Delta}^*$ is 
   $$\calD^{-\gamma}_{p^\prime}(\underline{\Delta}^*)=
     \calD^{-\gamma}_{p^\prime}(\Delta_{\min})\oplus
     \mathop{\mbox{$\oplus$}}_{q\in I_\gamma}\underline{\calE}_q^\perp,$$
  where the spaces $\underline{\calE}_{q_j^\pm}^\perp$ are defined as 
  follows: 
   \begin{itemize}
    \item[{\rm i)}] If either $q_j^\pm\not=0$ or 
     $\mbox{\rm dim}\,\bz\not=2$, there exists a unique subspace 
     $\underline{E}_j\subset E_j$ such that 
     $\underline{\calE}_{q_j^\pm}=\underline{E}_j\otimes\omega\,   
     t^{-q_j^\pm}$. Then we set 
      $$\underline{\calE}_{q_j^\pm}^\perp=\underline{E}_j^\perp\otimes
        \omega\,t^{-q_j^\mp},$$
     where $\underline{E}_j^\perp$ is the orthogonal complement of 
     $\underline{E}_j$ in $E_j$ with respect to the $L_2(\partial\bz)$-scalar
     product. 
    \item[{\rm ii)}] If $\mbox{\rm dim}\,\bz=2$ and $q_j^\pm=0$ define 
     $\underline{\calE}_0^\perp=\{0\}$ if $\underline{\calE}_0=\calE_0$, 
     $\underline{\calE}_0^\perp=\calE_0$ if $\underline{\calE}_0=\{0\}$, and 
     $\underline{\calE}_0^\perp=\underline{\calE}_0$ if 
     $\underline{\calE}_0=E_0\otimes\omega$. 
   \end{itemize}
  Note that $\underline{\calE}_{q_j^\pm}^\perp$ is a subspace of
  $\calE_{q_j^\mp}$ or, equivalently, $\underline{\calE}_{q}^\perp$ is a
  subspace of $\calE_{(n-1)-q}$.
 \end{theorem} 
  
 \begin{corollary}\label{selfadjoint}
  The selfadjoint dilation invariant extensions $\underline{\Delta}$ of
  $\Delta$ in $\calH^{0,0}_2(\bz)$ are those with a domain of the form 
   $$\calD^0_2(\underline{\Delta})=\calD^0_2(\Delta_{\min})\oplus
     \mathop{\mbox{$\oplus$}}_{q\in I_0}\underline{\calE}_q$$
  with $\underline{\calE}_{q}^\perp=\underline{\calE}_{(n-1)-q}$ for all 
  $q\in I_0$ $($in particular $\underline{\calE}_0=E_0\otimes\omega$ in case 
  $\mbox{\rm dim}\,\bz=2)$.
 \end{corollary}
 
 Applying Theorems 8.3 and 8.12 of \cite{GiMe}, the Friedrichs extension of
 $\Delta$ has the domain 
  $$\calD^0_2(\underline{\Delta})=
    \begin{cases}
     \calD^0_2(\Delta_{\min})\oplus
      \mathop{\mbox{$\oplus$}}\limits_{\substack{q\in I_0\\\re q<0}}
      \calE_q\oplus(E_0\otimes\omega) & 
      \mbox{\rm dim}\,\bz=2\\
     \calD^0_2(\Delta_{\min})\oplus
      \mathop{\mbox{$\oplus$}}\limits_{\substack
      {q\in I_0\\\re q\le\frac{n-1}{2}}}
      \calE_q & 
      \mbox{\rm dim}\,\bz\not=2  
    \end{cases}.$$ 
 In particular, the Friedrichs extension is dilation invariant. 

 \begin{remark}\label{gehtgenauso}
  All the results of Sections {\rm\ref{section4.2}} and {\rm\ref{section4.3}} 
  hold true in an analogous form for the model cone operator
  $\widehat{\Delta}$ considered as an unbounded operator in
  $\calK^{0,\gamma}_p(\partial\bz^\wedge)$. 
 \end{remark}
 
%%%%%%%%%%%%%%%%%%%%%%%%%%%%%%%%%%%%%%%%%%%%%%%%%%%%%%%%%%%%%%%%%%%%%

\subsection{Elliptic extensions}\label{section4.4}
 Proposition \ref{invdomain} provides a complete description of the closed
 extensions $\underline{\Delta}$ of $\Delta$ such that $-\underline{\Delta}$
 satisfies the ellipticity conditions {\rm(E1)} and {\rm(E2)} of Section
 \ref{section2.1}. We shall now discuss  condition {\rm(E3)}, assuming
that 
 $|\gamma|<\frac{1}{2}\mbox{\rm dim}\,\bz=\frac{n+1}{2}$
 (the choice of this range is
 connected to the scale of natural $L_p$-spaces on $\bz$ as we shall explain
 below). It turns out that for each given $\gamma$ we find at least one
 extension having property {\rm(E3)}; in case $\mbox{\rm dim}\,\bz\le3$ we
 find more than one. However, the extensions we describe might not represent
 all possible choices.

 \begin{theorem}\label{elldomain1} Let   $\Lambda\subset\cz\setminus\rz_+$
be an arbitrary sector. 
  Consider $-\Delta$ as an unbounded operator in $\calH^{0,\gamma}_p(\bz)$ and
  assume $\mbox{\rm dim}\,\bz\ge4$. Then conditions {\rm(E1)}, {\rm(E2)},
  and {\rm(E3)} of Section {\rm\ref{section2.1}}  are fulfilled
  by $-\Delta_{\max}$ in case
  $0\le\gamma<\frac{1}{2}\mbox{\rm dim}\,\bz$ and by $-\Delta_{\min}$ in case
  $-\frac{1}{2}\mbox{\rm dim}\,\bz<\gamma\le0$. 
 \end{theorem}
 
 The assumption on the dimension of $\bz$ in the previous theorem ensures that
 $\Delta$ in $\calH^{0,0}_2(\bz)$ is essentially self-adjoint or, in other
 words, the inverted conormal symbol has no pole in the interval $I_0$. We
 shall omit the proof of this theorem, since it is a simpler version of that
 for the following one (cf.\ also the proof of Theorem 7.1 in \cite{CSS1}).  
 
 \begin{theorem}\label{elldomain2}
  Consider $-\Delta$ as an unbounded operator in $\calH^{0,\gamma}_p(\bz)$, 
  assume $\mbox{\rm dim}\,\bz\le3$, and let 
  $|\gamma|<\frac{1}{2}\mbox{\rm dim}\,\bz$. An extension $-\underline{\Delta}$
  satisfies conditions {\rm(E1)}, {\rm(E2)}, and {\rm(E3)} of Section
  {\rm\ref{section2.1}} for any sector 
  $\Lambda\subset\cz\setminus\rz_+$, provided we choose its domain 
   $$\calD^\gamma_p(\underline{\Delta})=\calD^\gamma_p(\Delta_{\min})  
     \oplus\mathop{\mbox{$\oplus$}}_{q\in I_\gamma}\underline{\calE}_q$$
  according to the following rules: 
   \begin{itemize}
    \item[{\rm(i)}] If $q\in I_\gamma\cap I_{-\gamma}$, then 
     $\underline{\calE}_q^\perp=\underline{\calE}_{(n-1)-q}$. 
    \item[{\rm(ii)}] If $\gamma\ge0$ and $q\in I_\gamma\setminus I_{-\gamma}$,
     then $\underline{\calE}_q=\calE_q$. 
    \item[{\rm(iii)}] If $\gamma\le0$ and 
     $q\in I_{-\gamma}\setminus I_\gamma$, then $\underline{\calE}_q=\{0\}$. 
   \end{itemize}
  In particular, 
  $\calD^\gamma_p(\underline{\Delta})=\calD^\gamma_p(\Delta_{\max})$ if 
  $\gamma\ge1$ and 
  $\calD^\gamma_p(\underline{\Delta})=\calD^\gamma_p(\Delta_{\min})$ if 
  $\gamma\le-1$. 
 \end{theorem}
 
 \begin{proof}
  By Theorem \ref{independence2} we may assume that $p=2$, and by duality it 
  suffices to treat the case $\gamma\ge0$. Let $\underline{\Delta}_0$ denote
  the selfadjoint extension of $\Delta$ in $\calH^{0,0}_p(\bz)$ with
  $\calD^\gamma_2(\underline{\Delta})\subset 
  \calD^0_2(\underline{\Delta}_0)$. Such an extension always exists due to
  assumption {\rm(i)} on the domain of $\underline{\Delta}$ and by Corollary 
  \ref{selfadjoint}; its domain is 
   $$\calD^0_2(\underline{\Delta}_0)=\calD^0_2(\Delta_{\min})
     \oplus
     \mathop{\mbox{$\oplus$}}_{q\in I_0\setminus I_{-\gamma}}\calE_q
     \oplus
     \mathop{\mbox{$\oplus$}}_{q\in I_\gamma\cap I_{-\gamma}}
     \underline{\calE}_q.$$
  If we then pass to the associated model cone operators and use Remark
  \ref{gehtgenauso}, we get that 
   \begin{equation}\label{abcde}
    \lambda+\widehat{\underline{\Delta}}:
    \calD^\gamma_2(\widehat{\underline{\Delta}})\longrightarrow 
    \calK^{0,\gamma}_2(\partial\bz^\wedge),\qquad\lambda\not\in\rpbar,
   \end{equation}
  is injective, since 
  ${\rm spec}(-\widehat{\underline{\Delta}}_0)\subset\rpbar$ and 
  $\calD^\gamma_2(\widehat{\underline{\Delta}})\subset 
  \calD^0_2(\widehat{\underline{\Delta}}_0)$. 
  
  By Theorem \ref{adjoint} (in the formulation for model cone operators),
  the adjoint $\widehat{\underline{\Delta}}^*$ of 
  $\widehat{\underline{\Delta}}$ has the domain 
   $$\calD^{-\gamma}_2(\widehat{\underline{\Delta}}^*)=
     \calD^{-\gamma}_2(\widehat{\Delta}_{\min})\oplus 
     \mathop{\mbox{$\oplus$}}_{q\in I_\gamma\cap I_{-\gamma}}
     \underline{\calE}_q.$$ 
  Now let $\lambda\in\cz\setminus\rpbar$ and  
  $u\in\calD^{-\gamma}_2(\widehat{\underline{\Delta}}^*)$ be an element
  of the kernel of $\lambda+\widehat{\underline{\Delta}}^*$, i.e.\ 
  $(\lambda+\widehat{\Delta})u=0$. We shall show now that this implies $u=0$.
  To this end write $u=u_0+u_1$ with 
  $u_0\in\calD^{-\gamma}_2(\widehat{\Delta}_{\min})$ and 
  $u_1\in\mathop{\mbox{$\oplus$}}\limits_{q\in I_\gamma\cap I_{-\gamma}}
  \underline{\calE}_q$. Note that $u_0,u_1\in\calK^{0,0}_2(\partial\bz^\wedge)$
  by the assumption on the dimension of $\bz$ and the structure of the domain
  of $\underline{\Delta}$. Since $\widehat{\Delta}u_1\in\cicomp(\intb)$ (as
  this is true for any linear combination of functions from the spaces 
  $\calE_q$), we
  obtain $(\lambda+\widehat{\Delta})u_0=-
  \lambda u_1-\widehat{\Delta} u_1\in
  \calK^{0,0}_2(\partial\bz^\wedge)$. But this means that 
   $$u_0\in\calD^0_2(\widehat{\Delta}_{\max})\cap
     \calD^{-\gamma}_2(\widehat{\Delta}_{\min})=
     \calD^0_2(\widehat{\Delta}_{\min})\oplus 
     \mathop{\mbox{$\oplus$}}_{q\in I_0\setminus I_{-\gamma}}\calE_q
     \subset\calD^0_2(\widehat{\underline{\Delta}}_0).$$
  The last inclusion is valid by construction of $\underline{\Delta}_0$. This
  yields $u\in\calD^0_2(\widehat{\underline{\Delta}}_0)$ and
  $(\lambda+\widehat\Delta)u=0$, hence $u=0$, since 
  ${\rm spec}(-\widehat{\underline{\Delta}}_0)\subset\rpbar$. 
  
  This shows the bijectivity of \eqref{abcde}, since there 
  $\lambda+\widehat{\underline{\Delta}}$ is a Fredholm operator (this
  follows from  \cite{Lesc}, Proposition 1.3.16, together with a parametrix
  construction on $\partial\bz^\wedge$ as in the proof of Proposition
  \ref{elmar}), hence has closed range. 
 \end{proof}

%%%%%%%%%%%%%%%%%%%%%%%%%%%%%%%%%%%%%%%%%%%%%%%%%%%%%%%%%%%%%%%%%%%%%

\subsection{The Cauchy Problem}\label{section4.5}
 
 Let $1<p<\infty$ and let $L_p(\bz)$ denote the $L_p$-space on $\intb$
 associated to the measure induced by the conical metric on $\intb$. Then 
  $$L_p(\bz)=\calH^{0,\gamma_p}_p(\bz),\qquad
    \gamma_p=(n+1)\big(\mbox{$\frac{1}{2}-\frac{1}{p}$}\big).$$
 In fact, away from the boundary these spaces coincide by definition; thus it
 suffices to consider functions supported close to the boundary. But then, cf.\
 Definition \ref{sobolev}, 
  $$\|u\|_{\calH^{0,\gamma_p}_p(\bz)}^p=
    \int_{[0,1]\times\partial\bz}|t^{\frac{n+1}{2}-\gamma_p}u(t,x)|^p\,
    \mbox{$\frac{dt}{t}dx$}=
    \int_{[0,1]\times\partial\bz}|u(t,x)|^p\,t^ndtdx=
    \|u\|_{L_p(\bz)}^p.$$
 Clearly, $|\gamma_p|<\frac{n+1}{2}=\frac{1}{2}\mbox{\rm dim}\,\bz$ when $p$
 ranges from $1$ to $\infty$. Therefore the results of the previous Section
 \ref{section4.4} can be applied to the Laplacian in $L_p(\bz)$, $1<p<\infty$. 
 
 Combining these results with Theorem \ref{powers} and the Dore-Venni theorem
 (Theorem 3.2 in \cite{DoVe}), one obtains maximal regularity for solutions of
 the Cauchy problem: 
 
 \begin{theorem}
  Consider $\Delta$ as an unbounded operator in $L_p(\bz)$, $1<p<\infty$. For
  each closed extension $\underline{\Delta}$ from Theorems
  {\rm\ref{elldomain1}} or {\rm\ref{elldomain2}} associated with
  $\gamma=\gamma_p$, the Cauchy problem 
   $$\dot{u}(t)-\Delta u(t)=f(t)\quad\text{on}\quad 0<t<T,\qquad u(0)=0,$$
  has for any $f\in L_q([0,T],L_p(\bz))$, $1<q<\infty$, a unique solution 
   $$u\,\in\,W^1_q([0,T],L_p(\bz))\,\cap\,
     L_q([0,T],\calD^{\gamma_p}_p(\underline{\Delta})).$$ 
 \end{theorem}

%%%%%%%%%%%%%%%%%%%%%%%%%%%%%%%%%%%%%%%%%%%%%%%%%%%%%%%%%%%%%%%%%%%%%
%%%%%%%%%%%%%%%%%%%%%%%%%%%%%%%%%%%%%%%%%%%%%%%%%%%%%%%%%%%%%%%%%%%%%

\section{Appendix: Parameter-dependent cone pseudodifferential operators}
 \label{appendix}

 We try to give a concise review of the calculus of parameter-dependent 
 pseudodifferential operators on $\bz$ introduced by Schulze
 \cite{Schu1}, \cite{Schu2}. Our presentation follows \cite{Seil1} and
 \cite{GSS}. While there the parameter-space was $\rz^q$, we focus here on a
 subsector of the complex plane. The proofs pass over to this situation without
 any changes, and thus will be dropped here.
 
 We split the presentation into two parts: In Sections \ref{appendix1} to
 \ref{appendix3} we describe a sub-calculus of {\em flat} operators. 
 Under suitable ellipticity assumptions it already allows the construction of a
 rough parametrix to $\eta^\mu-A$ for a $\mu$-th order cone differential
 operator $A$. To describe the resolvent $(\eta^\mu-A)^{-1}$ we need to enlarge
 this calculus. This shall be explained starting with Section \ref{appendix4}. 
 
 In the following, $\Sigma$ is a closed sector in the complex plane
 (identified with $\rz^2$) containing zero, i.e.\ 
  $$\Sigma=\{\eta\in\cz\st
    \theta_1\le\text{\rm arg}\,\eta\le\theta_2\}\cup\{0\},\qquad
    -\pi\le\theta_1,\theta_2\le\pi.$$
 
 For a Fr\'echet space $E$, we let $\ci(\Sigma,E)$ denote the space of all
 continuous functions $\Sigma\to E$ that are smooth in the interior of $\Sigma$
 and whose derivatives have continuous extensions to the whole sector $\Sigma$.
 A subspace is $\calS(\Sigma,E)$, consisting of those functions that 
 decay rapidly in $\eta$ as $|\eta|\to\infty$ in $\Sigma$. 
 
 We fix a smooth positive function $[\cdot]$ with $[\eta]=|\eta|$ if
 $|\eta|\ge1$. Also recall that a cut-off function is a non-negative decreasing
 function in $\cicomp([0,1[)$ which is identically 1 near zero. 

%%%%%%%%%%%%%%%%%%%%%%%%%%%%%%%%%%%%%%%%%%%%%%%%%%%%%%%%%%%%%%%%%%%%%%

\subsection{Smoothing elements of the flat calculus}\label{appendix1}

 The space $\calC^{\infty,\infty}(\bz)$, consisting of all functions that are
 smooth in the interior of $\bz$ and vanish to infinite order at the boundary,
 is Fr\'echet in a natural way. Taking the projective tensor product yields the
 space 
  $$\calC^{\infty,\infty}(\bzbz)=\calC^{\infty,\infty}(\bz)\pit
    \calC^{\infty,\infty}(\bz).$$
 
 \begin{definition}\label{flatglobal}
  Let $C^{-\infty}_G(\Sigma)_\infty$ be the space of all operator-families 
  $r(\eta):\cii(\bz)\to\cii(\bz)$, $\eta\in\Sigma$, such that 
   $$(r(\eta)u)(y)=\int_\bz k_r(\eta,y,y^\prime)u(y^\prime)\,dy^\prime,$$
  where $dy'$ is a measure induced by a conic metric on $\bz$ and the kernel
  $k_r\in\calS(\Sigma,\cii(\bzbz))$ is rapidly decreasing in
  $\eta\in\Sigma$. 
 \end{definition}
 
 Besides this kind of smoothing operators -- which act globally on $\bz$ and
 decay rapidly in the parameter -- we shall also need a class of
 smoothing operators that are localized near the boundary but have a
 non-trivial dependence on $\eta\in\Sigma$. 
 
 To this end let $\calS^\infty(\partial\bz^\wedge)$ denote the space of smooth
 functions $\rz_+\times\partial\bz\to\cz$ that vanish to infinite order in
 $t=0$ and decrease rapidly for $t\to\infty$. We then define 
  $$\calS^\infty(\partial\bz^\wedge\times\partial\bz^\wedge)=
    \calS^\infty(\partial\bz^\wedge)\pit\calS^\infty(\partial\bz^\wedge).$$
 
 \begin{definition}\label{flatgreen}
  Let $R^\mu_G(\Sigma)_\infty$, $\mu\in\rz$, denote the space of all
  operator-families $a(\eta):\calS^\infty(\partial\bz^\wedge)\to
  \calS^\infty(\partial\bz^\wedge)$, $\eta\in\Sigma$, such that 
   \begin{equation}\label{kernel}
     (a(\eta)u)(t,x)=[\eta]^{n+1}\int_{\partial\bz^\wedge}
     k_a(\eta,t[\eta],x,t^\prime[\eta],x^\prime)u(t^\prime,x^\prime)\,
     {t^\prime}^{n}dt^\prime dx^\prime,
   \end{equation}
  with an integral kernel satisfying 
   $$k_a(\eta,t,x,t^\prime,x^\prime)\in
     S^\mu_{cl}(\Sigma,
     \calS^\infty(\partial\bz^\wedge\times\partial\bz^\wedge)):=
     S^\mu_{cl}(\Sigma)\pit
     \calS^\infty(\partial\bz^\wedge\times\partial\bz^\wedge).$$
 \end{definition}
 
 Using such operator-families, the so-called 
 {\em flat Green symbols} 
%% or parameter-dependent {\em flat Green operators} 
 are defined as follows: 

 \begin{definition}\label{flatgreen2}
  For $\mu\in\rz$ let $C^\mu_G(\Sigma)_\infty$ denote the space of all
  operator-families 
  $g(\eta):\calC^{\infty,\infty}(\bz)\to\calC^{\infty,\infty}(\bz)$, 
  $\eta\in\Sigma$, such that 
   \begin{equation}\label{flatgreensymb}
    g(\eta)=\sigma\,a(\eta)\,\sigma_0+r(\eta)
   \end{equation}
  for some cut-off functions $\sigma,\sigma_0\in\ci([0,1[)$, 
  $a\in R^\mu_G(\Sigma)_\infty$, and $r\in C^{-\infty}_G(\Sigma)_\infty$. 
 \end{definition}
 
 Note that if $g$ is as in \eqref{flatgreensymb}, then 
 $g(\eta)=\widetilde{\sigma}\,a(\eta)\,\widetilde{\sigma}_0+
 \widetilde{r}(\eta)$ for any
 choice of cut-off functions 
 $\widetilde{\sigma},\widetilde{\sigma}_0\in\ci([0,1[)$
 with a resulting $\widetilde{r}\in C^{-\infty}_G(\Sigma)_\infty$. Moreover, 
 the (pointwise) composition of such operator-families yields a map 
  \begin{equation}\label{greenalg}
   C^{\mu_0}_G(\Sigma)_\infty\times C^{\mu_1}_G(\Sigma)_\infty
   \longrightarrow C^{\mu_0+\mu_1}_G(\Sigma)_\infty.
  \end{equation}

%%%%%%%%%%%%%%%%%%%%%%%%%%%%%%%%%%%%%%%%%%%%%%%%%%%%%%%%%%%%%%%%%%%%%%

\subsection{Holomorphic Mellin symbols}\label{appendix2}

 A holomorphic Mellin symbol of order $\mu\in\rz$ is a function
 $h:\rpbar\times\cz\to L^\mu_{cl}(\partial\bz;\Sigma)$ depending smoothly on
 $t\in\rpbar$ and holomorphically on $z\in\cz$. It has its values in the
 Fr\'echet space of parameter-dependent pseudodifferential operators on the
 boundary of $\bz$. Moreover we require that 
  $$c_l(\delta):=\sup_{t\ge 0}[t]^l\trinorm{\partial_t^l
    h(t,\delta+i\tau)}$$
 is a locally bounded function of $\delta\in\rz$ for any $l\in\nz_0$ and any 
 semi-norm $\trinorm{\cdot}$ of $L^\mu_{cl}(\partial\bz;\rz_\tau\times\Sigma)$.
 We denote the space of all such symbols by
 $M^\mu_O(\rpbar\times\partial\bz;\Sigma)$ and write 
$M^\mu_O(\partial\bz;\Sigma)$ for the subspace of  $t$-independent
 symbols.

 With $h\in M^\mu_O(\rpbar\times\partial\bz;\Sigma)$ we associate an
 operator-family 
 $\calS^\infty(\partial\bz^\wedge)\to\calS^\infty(\partial\bz^\wedge)$ by 
  \begin{equation}\label{mellinop}
   (\op_M(h)(\eta)u)(t,x)=\int_\Gamma t^{-z}h(t,z,t\eta)(\calM u)(z,x)
   \,\dbar z,
   \qquad u\in\calS^\infty(\partial\bz^\wedge),
  \end{equation}
 where $\Gamma$ is an arbitrary vertical line in the complex plane (the 
  arbitrariness is due to the holomorphy of $\calM u$ and Cauchy's
integral formula). Note that on the
 right-hand side of \eqref{mellinop} we do not use the symbol $h(t,z,\eta)$
 itself, but the `degenerate' one $h(t,z,t\eta)$. We refer to operators of that
 kind as parameter-dependent Mellin pseudodifferential operators
 or, shortly, Mellin operators.
  
 \begin{remark}\label{mellinsymbofa}
  If $A$ is a cone differential operator as in \eqref{coneoperator} then, for
  any
  $\varphi\in\ci([0,1[)$,
   $$\varphi(\eta^\mu-A)=\varphi\,t^{-\mu}\,\op_M(h)(\eta),\qquad
     h(t,z,\eta)=\eta^\mu-\smsum_{j=0}^\mu a_j(t)z^j.$$
 \end{remark}
 
 Mellin operators behave well under composition: 
 
 \begin{theorem}
  Let $h_j\in M^{\mu_j}_O(\rpbar\times\partial\bz;\Sigma)$ for $j=0,1$. Then  
  \begin{equation}\label{leibniz}
   (h_0\#h_1)(t,z,\eta)=\iint s^{i\tau}
                        h_0(t,z+i\tau,\eta)h_1(st,z,s\eta)
                        \,\frac{ds}{s}\dbar\tau
  \end{equation}
  defines an element 
  $h_0\#h_1\in M^{\mu_0+\mu_1}_O(\rpbar\times\partial\bz;\Sigma)$, 
  the so-called Leibniz product, and 
   $$\op_M(h_0)(\eta)\op_M(h_1)(\eta)=\op_M(h_0\#h_1)(\eta)\qquad
     \forall\;\eta\in\Sigma.$$
 \end{theorem}
 
 The right-hand side of \eqref{leibniz} is understood as an oscillatory
 integral in a suitable sense. 

%%%%%%%%%%%%%%%%%%%%%%%%%%%%%%%%%%%%%%%%%%%%%%%%%%%%%%%%%%%%%%%%%%%%%%

\subsection{The calculus of flat cone operators}\label{appendix3}
 
 The operator-families we now consider are, roughly speaking, those which 
 are usual parameter-dependent pseudodifferential operators in the interior of
 $\bz$, and which are parameter-dependent Mellin operators near the boundary. 
 The global smoothing elements are flat Green symbols. Let us make this
 precise:
 
 \begin{definition}\label{flatcone}
  Let $\mu\in\rz$. Then $C_O^\mu(\Sigma)$ denotes the space of all
  operator-families $\cii(\bz)\to\cii(\bz)$ of the form 
  \begin{equation}\label{flatconesymb}
   c(\eta)=\sigma\,t^{-\mu}\,\op_M(h)(\eta)\,\sigma_0\;+\;
           (1-\sigma)\,p(\eta)\,(1-\sigma_1)\;+\;g(\eta),
  \end{equation}
  where $\sigma$, $\sigma_0$, $\sigma_1$ are cut-off functions
  satisfying $\sigma\sigma_0=\sigma$, $\sigma\sigma_1=\sigma_1$, and 
  \begin{itemize}
   \item[a)] $h(t,z,\eta)\in M^\mu_O(\rpbar\times\partial\bz;\Sigma)$
    is a holomorphic Mellin symbol, cf.\ Section \mbox{\rm\ref{appendix2}},  
   \item[b)] $p(\eta)\in L^\mu_{cl}(2\bz;\Sigma)$ is a parameter-dependent
    pseudodifferential operator on $2\bz$, 
   \item[c)] $g(\eta)\in C^\mu_G(\Sigma)_\infty$ is a Green symbol, cf.\
    Definition {\rm \ref{flatgreen}}.
  \end{itemize}
 \end{definition}

 For any choice of $0<\varrho<1$ one can achieve that  
 the symbols $h$ and $p$ in the representation \eqref{flatconesymb} are 
 {\em compatible} in the sense that 
  $$\varphi\left\{t^{-\mu}\,\op_M(h)(\eta)-p(\eta)\right\}\psi\;\in\;
    C^{-\infty}_G(\Sigma)_\infty\qquad
    \forall\;\varphi,\psi\in\cicomp(]\varrho,1[).$$
 In order to formulate the calculus in a smooth way, we shall fix such a
 $\varrho$ and shall always assume this compatibility relation to be satisfied.
 Moreover, we assume the involved cut-off functions $\sigma,\sigma_0,\sigma_1$
 to be identically 1 in a neighborhood of $[0,\varrho]$. Occasionally, we
 shall write $c(\eta)=\op(h,p,g)$ if $c(\eta)$ is as in $\eqref{flatconesymb}$.

 \begin{theorem}\label{flatalgebra}
  The pointwise composition of operator-families yields a map 
   $$C_O^{\mu_0}(\Sigma)\times C_O^{\mu_1}(\Sigma)
     \longrightarrow C_O^{\mu_0+\mu_1}(\Sigma).$$
  More precisely, if $c_j(\eta)=\op(h_j,p_j,g_j)$ for $j=0,1$, then
   $$c_0(\eta)c_1(\eta)=\op((T^{\mu_1}h_0)\#h_1,p_0p_1,\widetilde{g})$$ 
  with a resulting Green symbol
  $\widetilde{g}\in C_G^{\mu_0+\mu_1}(\Sigma)_\infty$. Recall
that $(T^\delta h)(t,z,\eta)=h(t,z+\delta,\eta)$. 
 \end{theorem}

 The operator-families from $C_O^{\mu}(\Sigma)$ introduced above are a
 subclass of parameter-dependent pseudodifferential operators on the
 interior of $\bz$. In particular, we can associate with them the usual
 homogeneous principal symbol 
  \begin{equation}\label{principal}
   \sigma_\psi^\mu(c)(y,\varrho,\eta)\in
   \ci((T^*\intb\times\Sigma)\setminus0)
  \end{equation}
 with $(y,\varrho)$ referring to variables of the cotangent bundle of $\intb$.
 In the coordinates $y=(t,x)$ near the boundary with corresponding
 covariables $\varrho=(\tau,\xi)$, the principal symbol has the form 
  $$\sigma_\psi^\mu(c)(t,x,\tau,\xi,\eta)=
    t^{-\mu}\,a_{(\mu)}(t,x,t\tau,\xi,t\eta)$$
 with a function $a_{(\mu)}(t,x,\tau,\xi,\eta)$, which is smooth in
 $(t,x)\in\rpbar\times\rz^n$ and
 $0\not=(\tau,\xi,\eta)\in\rz^{n+1}\times\Sigma$, and is positive homogeneous
 of order $\mu$ in $(\tau,\xi,\eta)$. Passing to the symbol
 $a_{(\mu)}(0,x,\tau,\xi,\eta)$ globally leads to the definition of the 
 {\em rescaled} principal symbol 
  \begin{equation}\label{rescaled}
   \widetilde{\sigma}_\psi^\mu(c)(x,\tau,\xi,\eta)\in
   \ci((T^*\partial\bz\times\rz\times\Sigma)\setminus0).
  \end{equation}
 Roughly speaking, this rescaled symbol describes the behavior of the
 principal symbol in the conical singularity itself. 
 We say that $c$ is {\em $\bz$-elliptic} if 
 \begin{itemize}
  \item[(E)] both the principal symbol $\sigma_\psi^\mu(c)$ and the rescaled
   symbol $\widetilde{\sigma}_\psi^\mu(c)$ are pointwise everywhere invertible. 
 \end{itemize}
 This condition allows the construction of a rough parametrix: 
 
 \begin{theorem}\label{roughparametrix1}
  Let $c(\eta)=\op(h_0,p_0,g_0)\in C_O^{\mu}(\Sigma)$ be $\bz$-elliptic.
  Then there exists an operator-family $b(\eta)=\op(h_1,p_1,g_1)\in 
  C_O^{-\mu}(\Sigma)$ such that 
  \begin{align*}
   b(\eta)c(\eta)&=1+\omega(t[\eta])\,\op_M(f_L)(\eta)\,\omega_0(t[\eta])
                   +g_L(\eta)\\
   c(\eta)b(\eta)&=1+\omega(t[\eta])\,\op_M(f_R)(\eta)\,\omega_0(t[\eta])
                   +g_R(\eta)
  \end{align*}
  with an arbitrary choice of cut-off functions
  $\omega$, $\omega_0$, Mellin symbols 
  $f_L,f_R\in M^{-\infty}_O$ $(\rpbar\times\partial\bz;\Sigma)$, 
  and flat Green symbols $g_L,g_R\in C^0_G(\Sigma)_\infty$. Moreover,
  $f_L=(T^\mu h_1)\#h_0-1$ and $f_R=(T^{-\mu}h_0)\#h_1-1$ on $[0,1[$.
 \end{theorem}
 
 Hence, $\bz$-elliptic symbols can be inverted up to smoothing remainders. However,
 this parametrix is not quite satisfactory, since a smoothing Mellin term is
 present and the Green symbols still have order $0$. To improve the quality of
 the remainder, one has to enlarge the calculus substantially (and has to pose
 additional ellipticity conditions). The elements of this enlarged calculus
 will be described in the next sections. 

%%%%%%%%%%%%%%%%%%%%%%%%%%%%%%%%%%%%%%%%%%%%%%%%%%%%%%%%%%%%%%%%%%%%%%

\subsection{Green symbols with asymptotics}\label{appendix4}
 
 Let $E^0$, $E^1$ be Banach spaces and  
 $\kappa^j=\{\kappa^j_\varrho\st\varrho>0\}\subset\calL(E^j)$ a strongly
 continuous group on $E^j$, i.e.\ $\kappa^j_1=1$ and
 $\kappa^j_{\varrho}\kappa^j_{\sigma}=\kappa^j_{\varrho\sigma}$. We refer
 to $\kappa^j$ as the {\em group action} of $E^j$.

 A function $a$ in $\ci(\Sigma,\calL(E^0,E^1))$ is said to be a symbol of order
 $\mu\in\rz$, if 
  $$\|\kappa_{1/\spk{\eta}}^1\,
    \partial_\eta^\alpha a(\eta)\,\kappa^0_{\spk{\eta}}\|_{\calL(E^0,E^1)}\le 
    c_\alpha\,\spk{\eta}^{\mu-|\alpha|}$$
 uniformly in $\eta\in\Sigma$ and for all multi-indices $\alpha$. We 
 write $a\in S^\mu(\Sigma;E^0,E^1)$. 
 
 A function $a\in\ci(\Sigma\setminus\{0\},\calL(E^0,E^1))$ is called {\em
 $($twisted$)$ homogeneous} of order $\mu\in\rz$, if
  \begin{equation}\label{twisted}
   a(\varrho\eta)=\varrho^\mu\,\kappa_{\varrho}^1\,a(\eta)\,
   \kappa_{1/\varrho}^0\qquad \forall\;\varrho>0,\,\eta\not=0.
  \end{equation}
 We shall denote the space of such functions by $S^{(-\mu)}(\Sigma;E^0,E^1)$. 
 The standard concept of classical
 (polyhomogeneous) symbols  having asymptotic expansions into
 homogeneous components extends to this operator-valued situation, resulting in the space
 $S^\mu_{cl}(\Sigma;E^0,E^1)$. 
 
 As a straightforward modification, one can admit $E^1$ to be a Fr\'echet
 space, which is the projective limit of Banach spaces,
 $E^1=\varprojlim_{k\in\nz}E^1_k$ with  
 $E_1^1\hookleftarrow E_2^1\hookleftarrow\ldots$, such that the group action on
 $E^1_1$ induces (by restriction) the group actions on all $E^1_k$, $k\in\nz$.
 Then we simply set 
  $$S^\mu_{cl}(\Sigma;E^0,E^1)=\mathop{\mbox{\Large$\cap$}}_{k\in\nz}\,
    S^\mu_{cl}(\Sigma;E^0,E_k^1).$$
 In the sequel we shall introduce various distribution spaces on
 $\partial\bz^\wedge=\rz_+\times\partial\bz$. The group action $\kappa$ always
 will be that induced by 
  \begin{equation}\label{groupaction}
   (\kappa_\varrho u)(t,x)=\varrho^{\frac{n+1}{2}}u(\varrho t,x),\qquad
   u\in\cicomp(\partial\bz^\wedge).
  \end{equation}
 
 \begin{definition}
  Let $\gamma,\theta\in\rz$ and $\theta>0$. An {\em asymptotic type}
  $Q\in\mbox{\rm As}(\gamma,\theta)$
  is a finite set of triples $(q,l,L)$, where $q$ is a
  complex number with $\frac{n+1}{2}-\gamma-\theta<\re q<\frac{n+1}{2}-\gamma$,
  $l\in\nz_0$, and $L\subset\ci(\partial\bz)$ is a finite-dimensional space of
  smooth functions. We shall write $Q=O$ if $Q$ is the empty set. 
  
  The conjugate type $\overline{Q}\in\mbox{\rm As}(\gamma,\theta)$ to $Q$ 
is the set of triples $(\overline{q},l,L)$, where
  $(q,l,L)\in Q$. 
 \end{definition}
 
 With an asymptotic type $Q=\{(q_j,l_j,L_j)\st j=0,\ldots,N\}\in
 \mbox{\rm As}(\gamma,\theta)$ we associate a
 finite-dimensional subspace of smooth functions supported in
 $[0,1[\times\partial\bz$, namely 
  \begin{equation}\label{asymptotic}
   \calE_Q=\Big\{(t,x)\mapsto\omega(t)\smsum_{j=0}^N\smsum_{k=0}^{l_j}\,
   u_{jk}(x)\,t^{-q_j}\,\log^kt\st u_{jk}\in L_j\Big\}. 
  \end{equation}
 Here $\omega$ is an arbitrary cut-off function. 
 
 \begin{definition}\label{subspaces}
  Let $s\in\rz$ and $Q\in\mbox{\rm As}(\gamma,\theta)$ be an asymptotic
  type. Then define 
  \begin{align*} 
   \calK^{s,\gamma}_{p,Q}(\partial\bz^\wedge)&=\calE_Q 
     \oplus\varprojlim_{\varepsilon>0}
     \calK^{s,\gamma+\theta-\varepsilon}_p(\partial\bz^\wedge)\\
   \calS^\gamma_Q(\partial\bz^\wedge)&=
     \{u\in\calK^{\infty,\gamma}_{2,Q}(\partial\bz^\wedge)\st 
     (1-\omega)u\in\calS(\partial\bz^\wedge)\}.
  \end{align*}
 In case $Q$ is the empty set, we shall write 
 $\calK^{s,\gamma}_{p,\theta}(\partial\bz^\wedge)$ and 
 $\calS^\gamma_\theta(\partial\bz^\wedge)$, respectively. 
 \end{definition}

 The spaces are clearly independent of the choice of the cut-off
 function. Moreover, they are Fr\'echet and can be written as 
 projective limits of Banach spaces. 
 This allows us to introduce Green operators as operator-valued symbols in 
 the above sense. 
 
 \begin{definition}\label{edgegreen}
  Let $Q\in\mbox{\rm As}(-\gamma,\theta)$, 
  $Q^\prime\in\mbox{\rm As}(\gamma^\prime,\theta^\prime)$ be given asymptotic
  types. We denote by 
  $R^\mu_G(\Sigma;\gamma,\theta,Q;\gamma^\prime,\theta^\prime,Q^\prime)$
  the space of all functions 
  $a:\Sigma\to\calL(\calK^{0,\gamma}_2(\partial\bz^\wedge),
  \calK^{0,\gamma^\prime}_2(\partial\bz^\wedge))$ with 
   $$a\in\mathop{\mbox{\Large$\cap$}}_{s\in\rz}
       S^\mu_{cl}(\Sigma;\calK^{s,\gamma}_2(\partial\bz^\wedge),
       \calS^{\gamma^\prime}_{Q^\prime}(\partial\bz^\wedge)),\qquad
     a^*\in\mathop{\mbox{\Large$\cap$}}_{s\in\rz}
       S^\mu_{cl}(\Sigma;\calK^{s,-\gamma^\prime}_2
       (\partial\bz^\wedge),\calS^{-\gamma}_Q(\partial\bz^\wedge)).$$
  Here, the $*$ refers to the pointwise adjoint with respect to the scalar
  product of $\calK^{0,0}_2(\partial\bz^\wedge)$. Moreover, we set
   $$R^\mu_G(\Sigma;\gamma,\theta;\gamma^\prime,\theta^\prime)=
     \mathop{\mbox{\Large$\cup$}}_{Q,Q^\prime}
     R^\mu_G(\Sigma;\gamma,\theta,Q;\gamma^\prime,\theta^\prime,Q^\prime)$$
  and  write $R^\mu_G(\Sigma;\gamma,\gamma^\prime,\theta)$ if
  $\theta=\theta^\prime$. 
 \end{definition}

 As an example, the flat Green symbols in Definition \ref{flatgreen} are 
 symbols of that type, namely 
  $$R^\mu_G(\Sigma)_\infty=
    \mathop{\mbox{\Large$\cap$}}_{\gamma,\gamma^\prime,\theta,\theta^\prime}\,
    R^\mu_G(\Sigma;\gamma,\theta,O;\gamma^\prime,\theta^\prime,O).$$
 It is often important to know that Green symbols have integral kernels with a specific
 structure. Set 
  $$\calS^\gamma_0(\partial\bz^\wedge)=
     \{u\in\calK^{\infty,\gamma}_2(\partial\bz^\wedge)\st 
     (1-\omega)u\in\calS(\partial\bz^\wedge), 
     (\log^kt)\omega u\in\calK^{\infty,\gamma}_2(\partial\bz^\wedge)
     \;\forall\,k\in\nz_0\}.$$
 
 \begin{theorem}\label{characterization}
  Let $a:\Sigma\to\calL(\calK^{0,\gamma}_2(\partial\bz^\wedge),
  \calK^{0,\gamma^\prime}_2(\partial\bz^\wedge))$ for given asymptotic types  
  $Q\in\mbox{\rm As}(-\gamma,\theta)$ and  
  $Q^\prime\in\mbox{\rm As}(\gamma^\prime,\theta^\prime)$. Then 
  $a\in R^\mu_G(\Sigma;\gamma,\theta,Q;\gamma^\prime,\theta^\prime,Q^\prime)$
  if and only if $a$ satisfies \eqref{kernel} with a kernel 
     $$k_a\in S^\mu_{cl}(\Sigma)\pit
       \calS^{\gamma^\prime}_{Q^\prime}(\partial\bz^\wedge)
       \,\widehat{\otimes}_\Gamma\, 
       \calS^{-\gamma}_{\overline{Q}}(\partial\bz^\wedge),$$
    where we have set 
     $$\calS^{\gamma^\prime}_{Q^\prime}(\partial\bz^\wedge)
       \,\widehat{\otimes}_\Gamma\, 
       \calS^{-\gamma}_{\overline{Q}}(\partial\bz^\wedge)=
       [\calS^{\gamma^\prime}_{Q^\prime}(\partial\bz^\wedge)
       \pit\calS^{-\gamma}_0(\partial\bz^\wedge)]\cap
       [\calS^{\gamma^\prime}_0(\partial\bz^\wedge)\pit 
       \calS^{-\gamma}_{\overline{Q}}(\partial\bz^\wedge)].$$
 \end{theorem} 
 
 To define general Green symbols on $\bz$ we need to introduce some function
 spaces on $\bz$: 
  $$\calC^{\infty,\gamma}(\bz)=\{u\in\ci(\intb)\st 
      \omega u\in\calS^\gamma_0(\partial\bz^\wedge)\},\qquad
    \calC^{\infty,\gamma}_Q(\bz)=\{u\in\ci(\intb)\st 
      \omega u\in\calS^\gamma_Q(\partial\bz^\wedge)\}.$$
 These are subspaces of $\calH^{s,\gamma}_p(\bz)$, independent of the choice of
 the involved cut-off function $\omega$. We shall write
 $\calC^{\infty,\gamma}_\theta(\bz)$ if $Q=O\in\mbox{\rm As}(\gamma,\theta)$ is
 the empty asymptotic type. 
 
 Now we define 
 $C^{-\infty}_G(\Sigma;\gamma,\theta,Q;\gamma^\prime,\theta^\prime,Q^\prime)$
 as the space of all functions 
 $r:\Sigma\to\calL(\calH^{0,\gamma}_2(\bz),\calH^{0,\gamma^\prime}_2(\bz))$ 
 such that 
  $$r\in\mathop{\mbox{\Large$\cap$}}_{s\in\rz}
       S^{-\infty}(\Sigma;\calH^{s,\gamma}_2(\bz),
       \calC^{\infty,\gamma^\prime}_{Q^\prime}(\bz)),\qquad
    r^*\in\mathop{\mbox{\Large$\cap$}}_{s\in\rz}
       S^{-\infty}(\Sigma;\calH^{s,-\gamma^\prime}_2
       (\bz),\calC^{\infty,-\gamma}_Q(\bz)),$$
 where $*$ refers to the adjoint with respect to the scalar-product of
 $\calH^{0,0}_2(\bz)$ and all spaces are equipped with the trivial group action
 $\kappa\equiv1$. 
 
 Similar as in Definition \ref{flatglobal} above, such operator-families
 possess an integral kernel 
  $$k_r\in\calS(\Sigma)\pit
    \calC^{\infty,\gamma^\prime}_{Q^\prime}(\bz)
    \,\widehat{\otimes}_\Gamma\,
    \calC^{\infty,-\gamma}_{\overline{Q}}(\bz),$$ 
 where  
  $$\calC^{\infty,\gamma^\prime}_{Q^\prime}(\bz)
    \,\widehat{\otimes}_\Gamma\,
    \calC^{\infty,-\gamma}_{\overline{Q}}(\bz)=
    [\calC^{\infty,\gamma^\prime}_{Q^\prime}(\bz)\pit
     \calC^{\infty,-\gamma}(\bz)]\cap
    [\calC^{\infty,\gamma^\prime}(\bz)\pit
     \calC^{\infty,-\gamma}_{\overline{Q}}(\bz)].$$ 
 Taking the union over all possible asymptotic types leads to the spaces 
 $C^{-\infty}_G(\Sigma;\gamma,\theta;\gamma^\prime,\theta^\prime)$ and 
 $C^{-\infty}_G(\Sigma;\gamma,\gamma^\prime,\theta)$ if $\theta=\theta^\prime$.
 
 \begin{definition}\label{green}
  Let $C^\mu_G(\Sigma;\gamma,\theta;\gamma^\prime,\theta^\prime)$ denote the 
  space of all operator-families 
  $g(\eta):\calC^{\infty,\gamma}(\bz)\to\calC^{\infty,\gamma^\prime}(\bz)$, 
  $\eta\in\Sigma$, such that 
   $$g(\eta)=\sigma\,a(\eta)\,\sigma_0+r(\eta)$$
  for  cut-off functions $\sigma,\sigma_0,$ and
  $a\in R^\mu_G(\Sigma;\gamma,\theta;\gamma^\prime,\theta^\prime)$, 
  $r\in C^{-\infty}_G(\Sigma;\gamma,\theta;\gamma^\prime,\theta^\prime)$. 
 \end{definition}
 
 The pointwise composition of such operator-families  yields a map 
  $$C^{\mu_0}_G(\Sigma;\gamma^\prime,\theta_1^\prime;
      \gamma^{\prime\prime},\theta^{\prime\prime})
   \times C^{\mu_1}_G(\Sigma;\gamma,\theta;\gamma^\prime,\theta^\prime)
   \longrightarrow
   C^{\mu_0+\mu_1}_G(\Sigma;\gamma,\theta;
   \gamma^{\prime\prime},\theta^{\prime\prime}).$$

 Let us finish this subsection with a result we shall need for the proof of our
 main theorem. 

 \begin{lemma}\label{weight}
  Let $\gamma\in\rz$ and $0<\varepsilon<\frac{1}{2}$. Then 
   $$C^0_G(\Sigma;\gamma,1,O;\gamma,1,O)\subset 
     C^0_G(\Sigma;\gamma,1-2\varepsilon,O;
     \gamma+2\varepsilon,1-2\varepsilon,O).$$ 
 \end{lemma}
 \begin{proof}
  Let $g(\eta)\in C^0_G(\Sigma;\gamma,1,O;\gamma,1,O)$. By Definition
  \ref{green}, we can write $g(\eta)=\sigma\,a(\eta)\,\sigma_0+r(\eta)$ with  
  $a(\eta)\in R^0_G(\Sigma;\gamma,1,O;\gamma,1,O)$ and 
  $r(\eta)\in C^{-\infty}_G(\Sigma;\gamma,1,O;\gamma,1,O)$. We now have to show
  that 
   \begin{equation}\label{formel6}
    a(\eta)\in R^0_G(\Sigma;\gamma,1-2\varepsilon,O;
     \gamma+2\varepsilon,1-2\varepsilon,O) 
   \end{equation}
  and that $r(\eta)\in C^{-\infty}_G(\Sigma;\gamma,1-2\varepsilon,O;
  \gamma+2\varepsilon,1-2\varepsilon,O)$. We restrict ourselves to the proof of
  \eqref{formel6}, since the symbol $r(\eta)$ can be treated in an analogous,
  even simpler way. By Theorem \ref{characterization} it suffices to show that 
   \begin{align*}
    \calS^{\gamma}_1(\partial\bz^\wedge)\,
    \widehat{\otimes}_\Gamma\,
    \calS^{-\gamma}_1(\partial\bz^\wedge)&\subset
    \calS^{\gamma+2\varepsilon}_{1-2\varepsilon}(\partial\bz^\wedge)\,
    \widehat{\otimes}_\Gamma\,
    \calS^{-\gamma}_{1-2\varepsilon}(\partial\bz^\wedge)
% \\
%    &=[\calS^{\gamma+2\varepsilon}_{1-2\varepsilon}(\partial\bz^\wedge)\,
%    \widehat{\otimes}_\pi\,
%    \calS^{-\gamma}_0(\partial\bz^\wedge)]\,\cap\,
%    [\calS^{\gamma+2\varepsilon}_0(\partial\bz^\wedge)\,
%    \widehat{\otimes}_\pi\,
%    \calS^{-\gamma}_{1-2\varepsilon}(\partial\bz^\wedge)]
   \end{align*}
  (recall that we write $\calS^\gamma_\theta=\calS^\gamma_O$ if $O\in\mbox{\rm
  As}(\gamma,\theta)$ is the empty asymptotic type). Clearly, 
   $$\calS^{\gamma}_1(\partial\bz^\wedge)\,
     \widehat{\otimes}_\Gamma\,
     \calS^{-\gamma}_1(\partial\bz^\wedge)\subset
     \calS^{\gamma}_1(\partial\bz^\wedge)\,
     \widehat{\otimes}_\pi\,
     \calS^{-\gamma}_0(\partial\bz^\wedge)=
     \calS^{\gamma+2\varepsilon}_{1-2\varepsilon}(\partial\bz^\wedge)\,
     \widehat{\otimes}_\pi\,
     \calS^{-\gamma}_0(\partial\bz^\wedge),$$
  where the last identity follows directly from the definition of the involved
  spaces. By Proposition 4.5 of \cite{Seil2} (in the version for operators on
  $\partial\bz^\wedge$) we have 
   $$\calS^{\gamma}_1(\partial\bz^\wedge)\,
     \widehat{\otimes}_\Gamma\,
     \calS^{-\gamma}_1(\partial\bz^\wedge)=
     \mathop{\mbox{\Large$\cap$}}_{0\le\sigma\le1}\,
     \calS^{\gamma}_\sigma(\partial\bz^\wedge)\,
     \widehat{\otimes}_\pi\,
     \calS^{-\gamma}_{1-\sigma}(\partial\bz^\wedge).$$
 Inserting $\sigma=2\varepsilon+\delta$, $\delta>0$ small, yields 
   $$\calS^{\gamma}_1(\partial\bz^\wedge)\,
     \widehat{\otimes}_\Gamma\,
     \calS^{-\gamma}_1(\partial\bz^\wedge)\subset
     \calS^{\gamma}_{2\varepsilon+\delta}(\partial\bz^\wedge)\,
     \widehat{\otimes}_\pi\,
     \calS^{-\gamma}_{1-2\varepsilon-\delta}(\partial\bz^\wedge)\subset 
     \calS^{\gamma+2\varepsilon}_0(\partial\bz^\wedge)\,
     \widehat{\otimes}_\pi\,
     \calS^{-\gamma}_{1-2\varepsilon-\delta}(\partial\bz^\wedge).$$
  Passing to the intersection over all $\delta>0$ gives 
   $$\calS^{\gamma}_1(\partial\bz^\wedge)\,
     \widehat{\otimes}_\Gamma\,
     \calS^{-\gamma}_1(\partial\bz^\wedge)\subset
     \calS^{\gamma+2\varepsilon}_0(\partial\bz^\wedge)\,
     \widehat{\otimes}_\pi\,
     \calS^{-\gamma}_{1-2\varepsilon}(\partial\bz^\wedge).$$
  This finishes the proof. 
 \end{proof}

%%%%%%%%%%%%%%%%%%%%%%%%%%%%%%%%%%%%%%%%%%%%%%%%%%%%%%%%%%%%%%%%%%%%%%

\subsection{Meromorphic Mellin symbols}\label{appendix5}
 
 An {\em asymptotic type} for Mellin symbols $P$ is a set of triples $(p,n,N)$
 with $p\in\cz$, $n\in\nz_0$, and $N$ a finite-dimensional subspace of finite
 rank operators from $L^{-\infty}(\partial\bz)$. Moreover, we require that 
  $\pi_\cz P\cap\{z\in\cz\st-\delta\le\re z\le\delta\}$
 is a finite set for each $\delta>0$, where 
  $$\pi_\cz P=\{p\in\cz\st (p,n,N)\in P\text{ for some }n,\,N\}.$$
 We shall write $P=O$ if $P$ is the empty set. 

 A meromorphic Mellin symbol with asymptotic type $P$ is a meromorphic function
 $f:\cz\to L^{-\infty}(\partial\bz)$ with poles at most in the points of
 $\pi_\cz P$. Moreover it satisfies: If $(p,n,N)\in P$, then the principal part
 of the Laurent series of $f$ in $p$ is of the form
 $\sum\limits_{k=0}^nR_k(z-p)^{-k-1}$ with $R_k\in N$;
 if $\chi\in\ci(\cz)$ is a $\pi_\cz P$-excision function (i.e.\ identically
 zero in an $\varepsilon$-neighborhood around $\pi_\cz P$ and identically 1
 outside the $2\varepsilon$-neighborhood), then
 $c(\delta)=\trinorm{(\chi f)(\delta+i\tau)}$ is a locally bounded function in
 $\delta\in\rz$ for each semi-norm of
 $L^{-\infty}(\partial\bz;\rz_\tau)=\calS(\rz_\tau,L^{-\infty}(\partial\bz))$. 
 
 As in \eqref{mellinop} we can associate with meromorphic Mellin symbols a
 pseudodifferential operator. Now, however, the operator will depend on the choice
 of the line. Letting $\Gamma_{1/2-\delta}$ denote the
vertical line $\{\Re z = 1/2-\delta\}$ we
 define $\op_M^{\delta}(f)$ by 
  \begin{equation}\label{wxyz}
   (\op_M^{\delta}(f)u)(t,x)=\int_{\Gamma_{\frac{1}{2}-\delta}}t^{-z}
   f(t,z)(\calM u)(z,x)\,\dbar z.
  \end{equation}
 Of course, we  have to require that none of the poles of $f$ lies on the
 chosen line; this we  shall always assume implicitly.  
 
 \begin{definition}\label{mplusg}
  Let $\gamma,\mu\in\rz$ and $k\in\nz$. Then
  $C^\mu_{M+G}(\Sigma;\gamma,\gamma-\mu,k)$ denotes the space of all
  operator-families
  $\calC^{\infty,\gamma}(\bz)\to\calC^{\infty,\gamma-\mu}(\bz)$ of the form 
   \begin{equation}\label{mplusgsymb}
    \omega(t[\eta])\Big(\smsum_{j=0}^{k-1}\smsum_{|\alpha|=0}^j
    t^{-\mu+j}\op_M^{\gamma_{j\alpha}-\frac{n}{2}}(f_{j\alpha})
    \eta^\alpha\Big)\widetilde\omega(t[\eta])\;+\;g(\eta),
   \end{equation}
  where $\omega$, $\widetilde{\omega}$ are arbitrary cut-off
  functions, $g\in C^\mu_G(\Sigma;\gamma,\gamma-\mu,k)$, the $f_{j\alpha}$
are meromorphic Mellin symbols with 
%\in M^{-\infty}_{P_{j\alpha}}$ for 
certain asymptotic types $P_{j\alpha}$, and
%  weights
  $\gamma_{j\alpha}\in\rz$ with $\gamma-j\le\gamma_{j\alpha}\le\gamma$. 
 \end{definition}
 
 Changing the cut-off functions in 
 \eqref{mplusgsymb}
 only yields remainders in $C^\mu_G(\Sigma;\gamma,\gamma-\mu,k)$.

%%%%%%%%%%%%%%%%%%%%%%%%%%%%%%%%%%%%%%%%%%%%%%%%%%%%%%%%%%%%%%%%%%%%%%

\subsection{The calculus of cone pseudodifferential operators}\label{appendix6}

 For $\gamma,\mu\in\rz$ and $k\in\nz$ let 
  \begin{equation}\label{calculus}
   C^\mu(\Sigma;\gamma,\gamma-\mu,k)=
   C^\mu_O(\Sigma)\;+\;C^\mu_{M+G}(\Sigma;\gamma,\gamma-\mu,k)
  \end{equation} 
 with $C^\mu_O(\Sigma)$ from Definition \ref{flatcone} and  
 $C^\mu_{M+G}(\Sigma;\gamma,\gamma-\mu,k)$ as in 
Definition \ref{mplusg}. The 
 elements of that space are operator-families 
 $\calC^{\infty,\gamma}(\bz)\to\calC^{\infty,\gamma-\mu}(\bz)$. Pointwise
 composition induces a map 
  $$C^{\mu_0}(\Sigma;\gamma-\mu_1,\gamma-\mu_1-\mu_0,k)\times 
    C^{\mu_1}(\Sigma;\gamma,\gamma-\mu_1,k)\longrightarrow
    C^{\mu_0+\mu_1}(\Sigma;\gamma,\gamma-\mu_1-\mu_0,k).$$
 \begin{remark}\label{mapprop}
  Let $c(\eta)\in C^\mu(\Sigma;\gamma,\gamma-\mu,k)$. Then, for each fixed
  $\eta$, $c(\eta)$ induces continuous maps 
   $$c(\eta):\calH^{s,\gamma}_p(\bz)\longrightarrow
     \calH^{s-\mu,\gamma-\mu}_p(\bz)$$
  for any $s\in\rz$ and $1<p<\infty$. Moreover, to any asymptotic type
  $Q\in\mbox{\rm As}(\gamma,k)$ there exists a type 
  $Q^\prime\in\mbox{\rm As}(\gamma-\mu,k)$ such that 
   $$c(\eta):\calC^{\infty,\gamma}_Q(\bz)\longrightarrow
     \calC^{\infty,\gamma-\mu}_{Q^\prime}(\bz).$$
 \end{remark}
 
 The construction of the parameter-dependent parametrix in Section
 \ref{section2.1} implicitly relies on ideas from the edge calculus; in
 particular, in \eqref{formel4} the notion of principal edge symbol is used to
 obtain a description of $(\eta^\mu-\underline{\widehat{A}})^{-1}$. We shall
 therefore recall those symbolic structures. In order to put this into
 perspective, we shall give -- although this will not be needed in this paper
 -- two theorems that show how the different types of ellipticity imply the
 existence of parametrices of corresponding quality. 
 
 Let $c(\eta)\in C^\mu(\Sigma;\gamma,\gamma-\mu,k)$ be given; then 
  $$c(\eta)=\sigma\,t^{-\mu}\,\op_M^{\gamma-\frac{n}{2}}(h)(\eta)\,\sigma_0+
    (1-\sigma)\,p(\eta)\,(1-\sigma_1)+(m+g)(\eta),$$
 where the first two terms are as in \eqref{flatconesymb}
 and $(m+g)(\eta)$ is as in  \eqref{mplusgsymb}. Since $(m+g)(\eta)$ has, in
 particular, a smooth distributional kernel, $c(\eta)$ is a
 parameter-dependent pseudodifferential operator on the interior of $\bz$,
 and we can associate with it the principal symbol and rescaled symbol as in
 \eqref{principal} and \eqref{rescaled}, respectively. 
 
 The {\em principal edge symbol} is 
  \begin{equation}\label{edgesymbol}
   \sigma^\mu_\wedge(c)(\eta)=
   t^{-\mu}\,\op_M^{\gamma-\frac{n}{2}}(h_0)(\eta)+
   \omega(t|\eta|)\Big(\smsum_{j=0}^{k-1}\smsum_{|\alpha|=j}
    t^{-\mu+j}\op_M^{\gamma_{j\alpha}-\frac{n}{2}}(f_{j\alpha})
    \eta^\alpha\Big)\widetilde\omega(t|\eta|)+g_{(\mu)}(\eta),
  \end{equation} 
 where $h_0(z,\eta)=h(0,z,\eta)$ and $g_{(\mu)}(\eta)$ is the homogeneous
 principal symbol of $g(\eta)\in R^\mu_G(\Sigma;\gamma,\gamma-\mu,k)$. 
 We consider the principal edge symbol as an operator-family
\begin{equation}\label{principaledge}
\sigma^\mu_\wedge(c)(\eta):\calK^{s,\gamma}_p(\partial\bz^\wedge)
    \longrightarrow 
    \calK^{s-\mu,\gamma-\mu}_p(\partial\bz^\wedge),\qquad\eta\not=0,
  \end{equation} 
 for $s\in\rz$ and $1<p<\infty$. 

Finally, the {\em conormal symbol} of 
 $c(\eta)$ is the meromorphic function with values in  $L^\mu_{cl}(\partial\bz)$
given by   
  \begin{equation}\label{conormalsymb}
   \sigma^\mu_M(c)(z)=h(0,z,0)+f_{00}(z): 
   H^{s}_p(\partial\bz)\longrightarrow 
   H^{s-\mu}_p(\partial\bz),\qquad z\in\cz.
  \end{equation} 
% respectively a meromorphic function with values in 

 We shall call $c(\eta)\in C^\mu(\Sigma;\gamma,\gamma-\mu,k)$ {\em elliptic},
 if 
  \begin{itemize}
   \item[(E)] both $\sigma_\psi^\mu(c)$ and $\widetilde{\sigma}_\psi^\mu(c)$
    are pointwise everywhere invertible (i.e.\ $c(\eta)$ is $\bz$-elliptic), 
   \item[$(\mbox{\rm E}_\wedge)$] 
    the principal edge symbol $\sigma_\wedge^\mu(c)$ in \eqref{principaledge} is pointwise everywhere 
    invertible. 
  \end{itemize}
 Here, the second condition initially is required to hold for some $s$ and $p$;
 but then it holds for all. As soon as $c(\eta)$ satisfies  one of
 the conditions (E) or $(\mbox{\rm E}_\wedge)$, the  conormal symbol will be
 meromorphically invertible. It is bijective on the vertical line 
 $\Gamma_{\frac{n+1}{2}-\gamma}$ in case 
 $c(\eta)$ satisfies $(\mbox{\rm E}_\wedge)$. 
 
 \begin{theorem}\label{roughparametrix2}
  Assume $c(\eta)\in C^\mu(\Sigma;\gamma,\gamma-\mu,k)$ satisfies condition 
  \mbox{\rm (E)} and the conormal symbol is invertible on the line
  $\Gamma_{\frac{n+1}{2}-\gamma}$. Then there exists a 
  $b(\eta)\in C^{-\mu}(\Sigma;\gamma-\mu,\gamma,k)$ 
  such that 
   $$b(\eta)c(\eta)-1\;\in\;C^0_G(\Sigma;\gamma,\gamma,k),\qquad
     c(\eta)b(\eta)-1\;\in\;C^0_G(\Sigma;\gamma-\mu,\gamma-\mu,k).$$ 
  This $($still rough$)$ parametrix $b(\eta)$ is uniquely determined modulo 
  $C^{-\mu}_G(\Sigma;\gamma-\mu,\gamma,k)$. 
 \end{theorem}
 
 \begin{theorem}\label{parametrix}
  Let $c(\eta)\in C^\mu(\Sigma;\gamma,\gamma-\mu,k)$ be elliptic.  
  Then there exists a $b(\eta)\in C^{-\mu}(\Sigma;\gamma-\mu,\gamma,k)$ 
  such that 
   $$b(\eta)c(\eta)-1\;\in\;C^{-\infty}_G(\Sigma;\gamma,\gamma,k),\qquad
     c(\eta)b(\eta)-1\;\in\;     
     C^{-\infty}_G(\Sigma;\gamma-\mu,\gamma-\mu,k).$$ 
  The parametrix $b(\eta)$ is uniquely determined modulo 
  $C^{-\infty}_G(\Sigma;\gamma-\mu,\gamma,k)$. 
 \end{theorem}
 
%%%%%%%%%%%%%%%%%%%%%%%%%%%%%%%%%%%%%%%%%%%%%%%%%%%%%%%%%%%%%%%%%%%%%
%%%%%%%%%%%%%%%%%%%%%%%%%%%%%%%%%%%%%%%%%%%%%%%%%%%%%%%%%%%%%%%%%%%%%

\begin{small}
\bibliographystyle{amsalpha}

\begin{thebibliography}{99}

\bibitem{BrSe1}
J.\ Br\"uning and R.\ Seeley.
An index theorem for regular singular operators. 
{\em Amer.\ J.\ Math.} {\bf 110}:659-714, 1988.

\bibitem{Che}
J.\ Cheeger. 
On the spectral geometry of spaces with cone-like singularities. 
{\em Proc.\ Nat.\ Acad.\ Sci.\ USA} {\bf 76}:2103-2106, 1979.

\bibitem{CSS1}
 S.\ Coriasco, E.\ Schrohe, J.\ Seiler.
 Bounded imaginary powers of cone differential operators.
 {\em Math.\ Z.}, to appear. (Preprint math.AP/0106008).

\bibitem{CSS2}
 S.\ Coriasco, E.\ Schrohe, J.\ Seiler.
 Differential operators on conic manifolds: Maximal regularity and parabolic
 equations.
 {\em Bull.\ Soc.\ Roy.\ Sci.\ Li\`ege} {\bf 70}: 207-229, 2001.

\bibitem{DoVe}
 G.\ Dore, A.\ Venni. 
 On the closednes of the sum of two operators. 
 {\em Math.\ Z.} {\bf 196}: 189-201, 1987. 

\bibitem{EgSc}
 Yu.\ Egorov, B.-W.\ Schulze. 
 {\em Pseudodifferential Operators, Singularities, Applications.}, 
 Birkh\"auser Verlag 1997.

\bibitem{Gil}
 J.B.\ Gil. 
 Heat trace asymptotics for cone differential operators. PhD thesis,
 Potsdam 1998. 

\bibitem{GilMN}
 J.B.\ Gil. 
 Full asymptotic expansion of the heat trace for non-self-adjoint elliptic cone operators.
 To appear in  {\em  Math. Nachr.}

\bibitem{GilLoya}
 J.B.\ Gil, P.\ Loya.
 On the noncommutative residue and the heat trace expansion on conic manifolds
 To appear in  {\em Manuscripta Math.}

 
\bibitem{GiMe}
 J.B.\ Gil, G.A.\ Mendoza. 
 Adjoints of elliptic cone operators. 
 Preprint, Temple University, Philadelphia, 2001. To appear in 
 {\em Amer. J. Math.}

\bibitem{GSS}
 J.B.\ Gil, B.-W.\ Schulze, J.\ Seiler. 
 Cone pseudodifferential operators in the edge symbolic calculus. 
 {\em Osaka J.\ Math.} {\bf 37}: 221-260, 2000.
 
\bibitem{Lesc} 
 M.\ Lesch. 
 {\em Operators of Fuchs type, Conical Singularities, and Asymptotic Methods}. 
 Teubner-Texte Math.\ {\bf 136}, Teubner-Verlag, 1997. 

\bibitem{Loya1}
 P.\ Loya.
 Complex Powers of Differential Operators on Manifolds with 
 Conical Singularities. 
 To appear in {\em  J. Anal. Math}.

\bibitem{Loya2}
 P.\ Loya.
 On the Resolvent of Differential Operators on Conic Manifolds. 
 To appear in {\em Comm. Anal. Geom.}


\bibitem{RBM} 
 R.\ Melrose. 
 Transformation of boundary value problems. 
 {\em Acta Math.} {\bf 147}:149-236, 1981.

\bibitem{Moo}
 E.\ Mooers. Heat kernel asymptotics on manifolds with conic singularities.
{\em  J. Anal. Math.} {\bf 78}: 1-36, 1999. 

 \bibitem{Plam86}
   B.\ Plamenevskij. {\em Algebras of pseudodifferential operators}.
   Nauka, Moscow 1986 (in Russian).

\bibitem{Schr}
 E.\ Schrohe. 
 Spaces of weighted symbols and weighted Sobolev spaces on manifolds. 
 In H.O.\ Cordes, B.\ Gramsch, H.\ Widom (eds.),
 {\em Pseudodifferential operators}, Oberwolfach (1986), pp.\ 360-377.
 LN Math.\ {\bf 1256}, Springer Verlag, 1987. 

\bibitem{ScSc}
 E.\ Schrohe, B.W.\ Schulze. 
 Edge-degenerate boundary value problems on cones. 
 In {\em Evolution equations and their applications in physical and life
 sciences}, Bad Herrenalb (1998), pp.\ 159-173. 
 Marcel Dekker, 2001. 

\bibitem{ScSe} 
 E.\ Schrohe, J.\ Seiler. 
 Ellipticity and invertibility in the cone algebra on $L_p$-Sobolev spaces. 
 {\em Int.\ Eq.\ Oper.\ Th.} {\bf 41}: 93-114, 2001. 

\bibitem{Schu1}
 B.-W.\ Schulze. 
 The Mellin pseudo-differential calculus on manifolds with corners. 
 In H.\ Triebel et al.\ (eds.), 
 {\em Symposium `Analysis on Manifolds with Singularities', Breitenbrunn 
 1990}, 
 Teubner-Texte Math.\ {\bf 131}, Teubner-Verlag, 1992. 

\bibitem{Schu2}
 B.-W.\ Schulze. 
 {\em Pseudo-differential Operators on Manifolds with Singularities}.  
 North-Holland, 1991. 

\bibitem{Seel0}
 R.\ Seeley.
 Complex powers of an elliptic operator.
 In {\em Amer.\ Math.\ Soc.\ Proc.\ Symp.\ Pure Math.}, volume~{\bf 10}:
 288-307, 1967.
 
\bibitem{Seel1}
 R.\ Seeley. 
 The resolvent of an elliptic boundary problem. 
 {\em Amer.\ J.\ Math.} {\bf 91}: 889-920, 1969.

\bibitem{Seel2}
 R.\ Seeley. 
 Norms and domains of the complex powers $A_B^z$. 
 {\em Amer.\ J.\ Math.} {\bf 93}: 299-309, 1971.

\bibitem{Seil1}
 J.\ Seiler. 
 {\em Pseudodifferential Calculus on Manifolds with Non-compact Edges}. 
 PhD-thesis, Institut f\"ur Mathematik, Potsdam, 1997. 
 
\bibitem{Seil2}
 J.\ Seiler. 
 The cone algebra and a kernel characterization of Green operators. 
 In J.B.\ Gil et al.\ (eds.), 
 {\em Approaches to Singular Analysis}, pp.\ 1-29.
 Birkh\"auser Verlag, 2001. 

\end{thebibliography}

\end{small}

%%%%%%%%%%%%%%%%%%%%%%%%%%%%%%%%%%%%%%%%%%%%%%%%%%%%%%%%%%%%%%%%%%%%%
%%%%%%%%%%%%%%%%%%%%%%%%%%%%%%%%%%%%%%%%%%%%%%%%%%%%%%%%%%%%%%%%%%%%%

\end{document}